\documentclass[11pt,letterpaper]{amsart}

\pdfoutput=1

\usepackage{amsmath,amssymb,amsthm}
\usepackage{mathtools,url}
\usepackage{stmaryrd,xcolor,tikz-cd,textcomp}
\usepackage{enumitem}
\usepackage[top=1.25in, bottom=1in, left=1.25in, right=1.25in]{geometry}
\newcommand{\into}{\hookrightarrow}

\newcommand\GHilb{\operatorname{G-Hilb}}

\DeclareMathOperator{\Hom}{Hom}
\DeclareMathOperator{\Ext}{Ext}
\DeclareMathOperator{\End}{End}

\DeclareMathOperator{\Spec}{Spec}

\DeclareMathOperator{\img}{im}

\DeclareMathOperator{\Hilb}{Hilb}
\DeclareMathOperator{\Aut}{Aut}
\DeclareMathOperator{\ch}{ch}
\DeclareMathOperator{\orbch}{\wtilde{\ch}}
\DeclareMathOperator{\tr}{tr}
\DeclareMathOperator{\td}{td}
\DeclareMathOperator{\orbtd}{\wtilde{\td}}

\DeclareMathOperator{\Supp}{Supp}
\DeclareMathOperator{\Ind}{Ind}

\DeclareMathOperator{\characteristic}{char}
\DeclareMathOperator{\Quot}{Quot}
\DeclareMathAlphabet{\mathpzc}{OT1}{pzc}{m}{it}

\newcommand{\wtilde}{\widetilde}

\newcommand{\inprod}[1]{\langle#1\rangle}
\newcommand{\tb}{\textbf}
\newcommand{\orbv}{\wtilde{v}}
\newcommand{\orbw}{\wtilde{w}}

\newcommand{\orbP}{\wtilde{P}}
\newcommand{\orbp}{\wtilde{p}}

\newcommand{\D}{\mathrm{D^b}}
\newcommand{\rk}{\mathrm{rk}}
\newcommand{\Pic}{\mathrm{Pic}}

\newcommand{\Qcoh}{\mathrm{Qcoh}}
\newcommand{\Coh}{\mathrm{Coh}}

\newcommand{\Stab}{\mathrm{Stab}}

\newcommand{\GL}{\mathrm{GL}}
\newcommand{\PGL}{\mathrm{PGL}}

\newcommand{\reg}{\mathrm{reg}}

\newcommand{\Real}{\mathrm{Re}}

\newcommand{\Amp}{\mathrm{Amp}}
\newcommand{\ob}{\mathrm{ob}}

\newtheorem{thm}{Theorem}[section]
\newtheorem{cor}[thm]{Corollary}
\newtheorem{prop}[thm]{Proposition}
\newtheorem{lem}[thm]{Lemma}

\newenvironment{sketch}{
	\proof}{\endproof}

\theoremstyle{definition}
\newtheorem{defn}[thm]{Definition}
\newtheorem{defn-prop}[thm]{Definition/Proposition}
\newtheorem{ex}[thm]{Example}

\newtheorem{rmk}[thm]{Remark}
\newtheorem{notation}[thm]{Notation}
\newtheorem{qn}[thm]{Question}
\newtheorem{constrn}[thm]{Construction}
\newtheorem{situation}[thm]{Situation}

\newcommand{\CC}{\mathbb{C}}

\newcommand{\QQ}{\mathbb{Q}}
\newcommand{\RR}{\mathbb{R}}

\newcommand{\ZZ}{\mathbb{Z}}

\newcommand{\cA}{\mathcal{A}}

\newcommand{\cD}{\mathcal{D}}
\newcommand{\cE}{\mathcal{E}}
\newcommand{\cF}{\mathcal{F}}

\newcommand{\cH}{\mathcal{H}}

\newcommand{\cK}{\mathcal{K}}
\newcommand{\cL}{\mathcal{L}}
\newcommand{\cM}{\mathcal{M}}

\newcommand{\cO}{\mathcal{O}}
\newcommand{\cP}{\mathcal{P}}

\newcommand{\cT}{\mathcal{T}}

\newcommand{\cV}{\mathcal{V}}

\newcommand{\cX}{\mathcal{X}}

\newcommand{\tu}[1]{\textup{#1}}

\begin{document}
\title{Equivariant Moduli Theory on $ K3 $ Surfaces}
\author{Yuhang Chen}

\begin{abstract}
	In this paper we study equivariant moduli spaces of sheaves on a $ K3 $ surface $ X $ under a symplectic action of a finite group. We prove that under some mild conditions, equivariant moduli spaces of sheaves on $ X $ are irreducible symplectic manifolds deformation equivalent to Hilbert schemes of points on $ X $ via a connection between Gieseker and Bridgeland moduli spaces, as well as the derived McKay correspondence.
\end{abstract}

\maketitle
\tableofcontents

\section{Introduction}
\vspace{1pt}
Let $ X $ be a complex projective $ K3 $ surface. Its canonical bundle $ \omega_X $ is trivial and admits a nowhere-vanishing global section $ \theta, $ which is a holomorphic symplectic form such that at each point $ p $ of $ X $, there is a non-degenerate skew-symmetric pairing
\begin{equation*}
\theta_p = dx \wedge dy: T_pX \times T_pX \to \CC,
\end{equation*}
where $ T_pX $ is the complex tangent space at $ p $ and $ (x,y) $ are some local coordinates around $ p. $ The symplectic form $ \theta $ is unique up to a scalar. An automorphism $ \sigma $ of $ X $ is called symplectic if it preserves the symplectic form $ \theta $, i.e.,
$$ \sigma^*\theta = \theta, $$ 
or equivalently, for every point $ p $ of $ X, $
\begin{equation*}
\det(d\sigma_p) = 1,
\end{equation*}
which implies that if $ \sigma $ is of finite order $ n, $ then
\begin{equation*}
d\sigma_p (x,y) = (\lambda x, \lambda^{-1} y),
\end{equation*}
where $ \lambda $ is a primitive $ n $th root of unity.

It's known that (see for example, \cite[Section 15.1]{huybrechts2016lectures}) if a symplectic automorphism $ \sigma $ has finite order $ n, $ then $ n \leq 8, $ and in this case, the number of points fixed by $ \sigma $ is finite and only depends on $ n $ as shown in Table \ref{table_number_fixed_points}. These numbers were computed by Mukai via the Lefschetz fixed point formula and some identities in number theory in \cite{mukai1988finite}.
\begin{table}
	\centering
	\caption{The number of fixed points $ f_n $ of a symplectic automorphism on a complex projective $ K3 $ surface with finite order $ n $}\label{table_number_fixed_points}
	\begin{tabular}{c|l*{7}{c}}
		$ n $ & \tu{2} & \tu{3} & \tu{4} & \tu{5} & \tu{6} & \tu{7} & \tu{8} \\ \hline
		$ f_n $ & \tu{8} & \tu{6} & \tu{4} & \tu{4} & \tu{2} & \tu{3} & \tu{2}
	\end{tabular}
\end{table}
A symplectic automorphism $ \sigma $ of finite order determines a cyclic group which acts on $ X $ faithfully. In general, we can consider a finite group $ G $ acting on $ X $ faithfully and symplectically, which means $ G $ can be identified with a symplectic automorphism subgroup of $ X $, and hence the possible groups are limited by the intrinsic geometry of the $ K3 $ surface $ X $. 
All such possible symplectic actions were classified by Xiao in \cite{xiao1996galois}. For example, there are in total $ 80 $ possible groups (including the trivial one) with a maximal size of $ 960 $ that can occur as a subgroup of the symplectic automorphism group of a $ K3 $ surface. Therefore, symplectic actions on $ K3 $ surfaces are well understood by now.

On the other hand, moduli spaces of sheaves on $ K3 $ surfaces have been studied extensively by Beauville \cite{beauville1983varietes} , Huybrechts \cite{huybrechts1997birational}, Mukai \cite{mukai1984symplectic}, O'Grady \cite{o1997weight}, and Yoshioka \cite{yoshioka2001moduli} among others, 
partially motivated by physical theories such as $ N = 4 $ supersymmetric Yang-Mills theory and the related $ S $-duality \cite{vafa1994strong}. By imposing a suitable stability condition known as Gieseker stability via an ample line bundle $ H $ on $ X $ called a polarization, and choosing an element $ x $ in the numerical Grothendieck ring $ N(X) $, one can construct a projective scheme $ M(X,x) $ from a Quot scheme via the geometric invariant theory (GIT) such that $ M(X,x) $ is a good moduli space of semistable sheaves on $ X $ with numerical class $ x. $ Indeed, such a GIT construction for moduli spaces of sheaves exists for an arbitrary projective scheme. 

Unfortunately, $ M(X,x) $ is not a coarse moduli space since two semistable sheaves may share a point in $ M(X,x) $ when they are so-called $ S $-equivalent\footnote{Here $ S $ stands for Seshadri who first introduced the notion of $ S $-equivalence for vector bundles over a curve in \cite{seshadri1967space} where he called it strong equivalence.}. The GIT approach produces a quasi-projective scheme $ M^s(X,x) $ which is a coarse moduli space of stable sheaves with numerical class $ x. $ With a suitable choice of $ x $ and $ H, $ we have 
$$ M(X,x) = M^s(X,x), $$ 
i.e., there is no semistable sheaf with numerical class $ x $, and it follows that $ M(X,x) $ is deformation equivalent to a Hilbert scheme of points $ \Hilb^n(Y) $ on a $ K3 $ surface $ Y $, which is an irreducible symplectic manifold\footnote{A smooth complex projective variety $ X $ is an irreducible symplectic manifold if $ H^0(X, \Omega^2) = \CC \omega $ for a holomorphic symplectic form $ \omega $ and $ H^1(X, \cO_X) = 0. $}. Here $ Y $ can be taken to be the original $ K3 $ surface as $ \Hilb^n(X) $ is deformation equivalent to $ \Hilb^n(Y) $ for any $ K3 $ surface $ Y $. Therefore, the deformation type of a moduli space $ M(X,x) $ is determined by its dimension. For an element $ x $ in $ N(X) $, there is an associated Mukai vector $ v(x) = \ch(x) \sqrt{\td_X}$ in the numerical Chow ring $ R(X) $ of $ X $ which is equivalent to the Chern character $ \ch(x) $ and hence to the element $ x $ as well, but has the advantage that the intersection product on $ R(X) $ induces a Mukai pairing 
$$ \inprod{\cdot\ {,}\ \cdot}: R(X) \times R(X) \to \ZZ $$
such that we have
\begin{equation}\label{eq_dim}
\dim M(X,x) = 2 - \inprod{v(x)^2}
\end{equation}
as a consequence of the Hirzebruch-Riemann-Roch (HRR) formula. 


Now, it's natural to ask the following
\begin{qn}\label{qn_equivariant_moduli_space}
	What can we say about equivariant moduli spaces of sheaves on a $ K3 $ surface with symplectic automorphisms?
\end{qn}

Let $G$ be a finite subgroup of the symplectic automorphism group of $X.$ The category of $ G $-equivariant sheaves on $ X $ is equivalent to the category of sheaves on the quotient stack $ \cX = [X/G]. $ Therefore, we can consider moduli spaces of sheaves on $ \cX $ instead. By a GIT construction of Nironi in \cite{nironi2009moduli}, there exists a moduli space $ M^{(s)}(\cX,x) $ of semistable (resp. stable) sheaves on $ \cX $ with numerical class $ x $ in $ N(\cX) $. Roughly speaking, one chooses a polarization $ (\cH, \cV) $
on $ \cX $ to impose stability conditions, where $ \cH $ is a line bundle on $ \cX $ which descends to its coarse moduli space, and $ \cV $ is a generating sheaf on $ \cX $ which is a locally free sheaf on $ \cX $ such that it contains all irreducible representations of the automorphism group at each point on $ \cX. $ In our case, the geometric quotient $ X/G $ is a coarse moduli space of $ \cX $, so we can pullback an ample line bundle on $ X/G $ to a line bundle $ \cH $ on $ \cX $, and there is a canonical generating sheaf 
$$ \cV_\reg = \cO_\cX \otimes \rho_\reg $$ 
on $ \cX $ where $ \rho_\reg $ is the regular representation of $ G. $ Therefore, fixing an element $ x $ of $ N(\cX), $ we have a moduli space $ M^{(s)}(\cX,x) $ of semistable (resp. stable) sheaves on $ \cX $ similarly to the case of schemes.

The canonical morphism
$$ p: X \to \cX $$
is smooth, and hence induces an exact pullback functor
\begin{equation*}
p^*: \Coh(\cX) \to \Coh(X), \quad \cE = (E, \phi) \mapsto E.
\end{equation*}
The line bundle $ \cH $ on $ \cX $ pulls back to a $ G $-invariant line bundle $ H = p^*\cH $ on $ X $. The degree of a line bundle $ L $ on $ X $ is defined by 
$$ \deg(L) = \deg(l h) $$
in $ \ZZ, $ where $ l = c_1(L) $ and $ h = c_1(H) $ in $R^1(X) \cong \Pic(X),$ the Picard group of $X.$
Let $ N(\cX) $ and $ N(X) $ denote the numerical Grothendieck rings of $ \cX $ and $ X $ respectively. 
For a coherent sheaf $ \cE $ (resp. $ E $) on $ \cX $ (resp. $ X $), let $ \gamma(\cE) $ (resp. $ \gamma(E) $) denote its image in $ N(\cX) $ (resp. $ N(X) $).
Then there is a numerical pullback 
$$ p^N: N(\cX) \to N(X), \quad \gamma(\cE) = \gamma(E,\phi) \mapsto \gamma(E). $$ 
The action of $ G $ on $ X $ induces an action on $ N(X) $ and an action on $ R(X) $. Let $ N(X)^G $ and $ R(X)^G $ denote their $ G $-invariant subspaces respectively. 
An element $ x $ in $ N(\cX) $ determines an orbifold Mukai vector $ \orbv(x) $ in the complex numerical Chow ring $ R(I\cX)_\CC $ of the inertia stack $ I\cX. $
In Section \ref{sec_HRR_K3/G}, we derive an explicit orbifold HRR formula for $ \cX = [X/G] $ via an orbifold Mukai pairing 
$$ \inprod{\cdot\ {,}\ \cdot}_{I\cX}: R(I\cX)_\CC \times R(I\cX)_\CC \to \CC, $$ 
and use it to obtain an orbifold version of formula (\ref{eq_dim}). See \cite[Section 3.4]{chen2023orbifold} for the definitions of the orbifold Mukai vector and the orbifold Mukai pairing for a quotient Deligne-Mumford (DM) stack.

We will answer Question \ref{qn_equivariant_moduli_space} in the following
\begin{thm}\label{thm_main}
	Let $ x $ be an element in $ N(\cX) $ with $ y = p^N x $ in $ N(X)^G $. Denote the Mukai vector of $ y $ by $ v = (r,c_1,s) $ in $ R(X)^G $ with $ d = \deg(y). $ Suppose $ r > 0 $ and the following conditions are satisfied:
	\begin{enumerate}[font=\normalfont,leftmargin=2em]
		\item $ y $ is primitive and $ H $ is $ y $-generic.
		\item $ d > 0 $ or $ \gcd(r,d) = 1 $.
	\end{enumerate}
	Then $ M(\cX,x) = M^s(\cX,x). $ If $ M(\cX,x) $ is non-empty, then it is an irreducible symplectic manifold of dimension $ n = 2 - \inprod{\orbv(x)^2}_{I\cX} $ deformation equivalent to $ \Hilb^{n/2}(X) $.
\end{thm}

The proof of Theorem \ref{thm_main} makes use of three notions: 
\begin{enumerate}[font=\normalfont,leftmargin=3em]
	\item Bridgeland stability conditions
	\item Induced stability conditions
	\item Fourier-Mukai transform
\end{enumerate}
See Section \ref{sec_proof_main_thm} for a complete proof. We outline the proof here.
\begin{sketch}
	Let $ M_1 = M(\cX,x). $ Roughly speaking, the argument can be divided into three steps.
	
	\tb{Step 1.} Our starting point is that the $ G $-invariant ample line bundle $ H $ on the $ K3 $ surface $ X $ determines a real number $ t_0 $ such that for all $ t \geq t_0, $ there exist $ G $-invariant stability conditions $ \sigma_t $ in the distinguished component $ \Stab^\dagger(X) $ in the space of stability conditions on $ \D(X), $ and $ H $-stable sheaves in $ \Coh(X) $ are precisely $ \sigma_t $-stable complexes in $ \D(X). $ These $ G $-invariant stability conditions $ \sigma_t $ induce stability conditions $ \tau_t $ on the $ G $-equivariant derived category $ \D(X)_G \cong \D(\cX) $ by a technique developed in \cite{macri2009inducing}. We then show that $ \cH $-stable sheaves in $ \Coh(\cX) $ are precisely $ \tau_t $-stable complexes in $ \D(\cX). $ This makes it possible to identify the Gieseker moduli space $ M_1 $ with a Bridgeland moduli space $ M_2 $ of stable complexes in $ \D(\cX). $
	
	\tb{Step 2.} The $ G $-Hilbert scheme of free orbits on $ X $ gives the minimal resolution $ M $ of the surface $ X/G, $ and induces a Fourier-Mukai transform 
	$$ \Phi: \D(\cX) \xrightarrow{\sim} \D(M) $$
	such that under $ \Phi, $ the stability conditions $ \tau_t $ on $ \D(\cX) $ are mapped to stability conditions $ \Phi^S(\tau_t) $ on $ \D(M) $, and the Bridgeland moduli space $ M_2 $ on $ \cX $ is then identified with another Bridgeland moduli space $ M_3 $ on $ M. $
	
	\tb{Step 3.} By \cite[Proposition 6.1]{beckmann2020equivariant}, we know these stability conditions $ \Phi^S(\tau_t) $ are in the distinguished component $ \Stab^\dagger(M). $ Under $ \Phi, $ the primitive element $ x $ in $ N(\cX) $ is transformed to another primitive element $ \Phi^N(x) $ in $ N(M) $. In the same time we can choose $ t $ such that $ \Phi^S(\tau_t) $ is generic with respect to the element $ \Phi^N(x) $. 
	The main result of \cite{bottini2021stable} then says that the Bridgeland moduli space $ M_3 $ is deformation equivalent to a Hilbert scheme of points on a $ K3 $ surface, and hence so is our moduli space $ M_1. $
\end{sketch}

\subsection{Structure}
The paper is structured as follows.

In Section \ref{sec_moduli_sheaves_stacks}, we review Gieseker stability conditions for coherent sheaves on polarized projective schemes and stacks. We then focus on the case of a projective quotient stack $ [X/G] $ where $X$ is a projective scheme and $G$ is a finite group, and define $ G $-equivariant moduli spaces of stable sheaves on $ X. $

In Section \ref{sec_moduli_sheaves_K3/G}, we derive an explicit orbifold HRR formula for a quotient stack $ [K3/G] $ and use it to compute the dimensions of $ G $-equivariant moduli spaces of stable sheaves on a $ K3 $ surface $ X $. As a joyful digression, we apply the same formula to reproduce the number of fixed points on $ X $ when $ G \cong \ZZ/n\ZZ $ without using the Lefschetz fixed point formula. 
Next we review the derived McKay correspondence, the Bridgeland stability conditions on $ X, $ and the induced stability conditions on $ [X/G] $. We then prove Theorem \ref{thm_main} in details. 


In the appendix, we review the Mukai pairing and the HRR formula for proper smooth schemes.

\subsection{Conventions}

Throughout this paper, the base scheme $ S = \Spec \CC $ unless otherwise specified. For a ring $ R, $ let $ R_\CC = R \otimes_\ZZ \CC $ denote its extension of scalars from $ \ZZ $ to $ \CC. $ Let $ G $ be a group. The identity of $ G $ is denoted by $ 1. $ A representation of $ G $ means a linear representation of $ G, $ i.e., a finite-dimensional complex vector space $ V $ with a group homomorphism $ G \to \GL(V). $ 

\subsection{Acknowledgements}
I would like to thank my advisor Hsian-Hua Tseng for providing a free research environment and help when necessary. I want to thank Yunfeng Jiang for suggesting the problem of $ [K3/G] $ and constant discussions throughout the research process. I thank Promit Kundu for many discussions on stacks and moduli spaces of sheaves on $ K3 $ surfaces. I thank Yonghong Huang for discussions on stacks. I am also thankful to Hao Sun for discussions on Bridgeland stability conditions and checking Section \ref{sec_Bridgeland}. I am grateful to David Anderson and James Cogdell for useful comments on an early version of my paper. I wish to thank Johan de Jong, Dan Edidin, Daniel Huybrechts, Martijn Kool, Georg Oberdieck, Martin Olsson, and Alessandra Sarti for useful email correspondence. I also thank Nicolas Addington, Benjamin Call, Kathleen Clark, Deniz Genlik, Jean-Pierre Serre, Luke Wiljanen, and Yilong Zhang for pointing out various mistakes in earlier versions of my paper.
\vskip 1pt

\section{Equivariant moduli theory}\label{sec_moduli_sheaves_stacks}
\vspace{3pt}

\subsection{Moduli spaces of sheaves on projective schemes}\label{sec_moduli_sheaves_schemes}
In this section we review Gieseker stability of coherent sheaves on projective schemes over $ \CC $.

Let's first recall some notions which are used to define Hilbert polynomials of coherent sheaves on projective schemes.
\begin{defn}
	The dimension of a coherent sheaf $ E $ on a scheme $ X $ is the dimension of its support as a subscheme of $ X $, i.e.,
	\begin{equation*}
	\dim E = \dim \Supp E.
	\end{equation*}
\end{defn}
\begin{defn}[{\cite[Section I.7]{hartshorne1977algebraic}}]
	A polynomial $ P(z) $ in $ \QQ[z] $ is called \tb{numerical} if $ P(m) $ is in $ \ZZ $ for all $ m \gg 0 $ in $ \ZZ. $
\end{defn}
\begin{defn}
	Let $ X $ be a projective scheme. A choice of an ample line bundle $ H $ on $ X $ is called a \tb{polarization} of $ X. $ A pair $ (X, H) $ where $ H $ is an ample line bundle on $ X $ is called a \tb{polarized projective scheme}.
\end{defn}
Let $ (X, H) $ be a polarized projective scheme.
\begin{defn-prop}
	 Let $ E $ be a coherent sheaf on $ X. $ There is a unique numerical
	polynomial 
	$$ P_E(z) = a_n(E) z^n + \dots + a_1(E)z + a_0(E) $$
	in $ \QQ[z] $ such that
	\begin{equation*}
	P_E(m) = \chi(X, E\otimes H^{\otimes m})
	\end{equation*}
	in $ \ZZ $ for all $ m $ in $ \ZZ. $ We call $ P_E $ the \tb{Hilbert polynomial} of $ E $ with respect to $ H. $ If $ E = 0 $, then $ P_E = 0 $; else, $ n = \dim E $. 
	If $ E \neq 0, $ then the \tb{reduced (or normalized) Hilbert polynomial} of $ E $ is defined by the monic polynomial \footnote{It can also be defined by $ P_E(z)/n!a_n(E) $ as in \cite[Definition 1.2.3]{huybrechts2010geometry}. This will give the same notion of stability for sheaves.}
	\begin{equation*}
	p_E(z) = \frac{P_E(z)}{a_n(E)}
	\end{equation*}
	in $ \QQ[z] $.
\end{defn-prop}
\begin{notation}\label{notation_rk_deg}
	Let $ (X, H) $ be a connected smooth polarized projective scheme of dimension $ d $. There are two rings associated with $ X $: the numerical Grothendieck ring $ N(X) $ and the numerical Chow ring $ R(X). $ There is a \tb{rank map}, i.e., a ring homomorphism
	\begin{equation*}
	\rk: N(X) \to \ZZ
	\end{equation*}
	defined by
	\begin{equation*}
	\rk(x) = \rk(V) - \rk(W)
	\end{equation*}
	for any element $ x = \gamma(V) - \gamma(W) $ in $ N(X) \cong N^0(X).$
	There is also a \tb{degree map}, i.e., a group homomorphism 
	$$ \deg: N(X) \to \ZZ $$
	defined by
	\begin{equation*}
	\deg(x) = \deg(c_1(x) h^{d-1})
	\end{equation*}
	for all $ x $ in $ N(X), $ where $ h = c_1(H) $ in $ R^1(X) $. If $ E $ is torsion-free, then $ n = d $ and by the HRR theorem, we have 
	\begin{equation*}
	a_d(E) = \frac{\rk(E)\deg(H)}{d!}
	\end{equation*}
	in $ \QQ, $ where $ \deg(H) = \deg(h^d) $ in $ \ZZ. $
\end{notation}
\begin{rmk}
	The Hilbert polynomial of a sheaf in $ \Coh(X) $ is additive on short exact sequences in $ \Coh(X) $ since the Euler characteristic $ \chi(X, \ \cdot\ ) $ is. Therefore, the Hilbert polynomial of a coherent sheaf
	only depends on its numerical class in $ N(X) $, and hence descends to an additive map
	\begin{equation*}
	P: N(X) \to \QQ[z], \quad \gamma(E) \mapsto P_E.
	\end{equation*}
	When $ X $ is smooth, then the Hilbert polynomial $ P_E $ of a coherent sheaf $ E $ on $ X $ is determined by the Mukai vector $ v(E) $ in $ R(X)_\QQ $ since the map $ v: N(X) \to R(X)_\QQ$ is injective.
\end{rmk}
Now we can define stability of pure sheaves.
\begin{defn}
	Let $ E $ be a nonzero coherent sheaf on $ X. $ We say the sheaf $ E $ is \tb{pure} if $ \dim F = \dim E $ for all nonzero subsheaves $ F \subset E $ on $ X. $ 
\end{defn}
\begin{rmk}
	The support of a pure sheaf is pure dimensional, but the converse is not true. A pure sheaf is torsion-free on its support. Indeed, pure sheaves are a generalization of torsion-free sheaves.
\end{rmk}
We introduce a notation to compare polynomials in $ \QQ[z]. $
\begin{notation}
	Let $ p(z) $ and $ q(z) $ be two polynomials in $ \QQ[z]. $ The inequality $ p(z) \leq q(z) $ means $ p(z) \leq q(z) $ for $ z \gg 0. $
	Similarly, the strict inequality $ p(z) < q(z) $ means $ p(z) < q(z) $ for $ z \gg 0. $
\end{notation}
\begin{defn}
	Let $ E $ be a coherent sheaf on $ X. $ We say $ E $ is \tb{$ H $-(semi)stable}\footnote{The $ H $-stability was first introduced by Gieseker in \cite{gieseker1977moduli} for vector bundles on surfaces, and it is commonly called Gieseker $ H $-stability or Gieseker stability.} if $ E $ is a nonzero pure sheaf and
	\begin{equation*}
	p_F(z) \ (\leq)\footnote{There are two statements here: the inequality $ \leq $ is used to define $ H $-semistable sheaves, and the strict inequality $ < $ is used to define $ H $-stable sheaves.} \ p_E(z)
	\end{equation*}
	for all proper nonzero subsheaves $ F \subset E $ on $ X. $ When the ample line bundle $ H $ is understood, $ H $-semistable (resp. $ H $-stable) sheaves are also called \tb{semistable} (resp. \tb{stable}) sheaves. We say $ E $ is \tb{strictly semistable} if it is semistable but not stable; in this case, there is a proper subsheaf $ F $ of $ E $ with $ p_F = p_E, $ and we say $ F $ is a destablizing subsheaf of $ E $ or destabilizes $ E. $
\end{defn}

Fix a numerical polynomial $ P $ in $ \QQ[z] $. Via a GIT construction, there is a projective scheme $ M_H(X,P) $ parametrizing $ S $-equivalent classes of semistable sheaves on $ X $ with Hilbert polynomial $ P $, which compactifies a quasi-projective scheme $ M_H^s(X,P) $ parametrizing stable sheaves on $ X $ with Hilbert polynomial $ P $. 

We can also fix a numerical class instead of a Hilbert polynomial. Let $ x $ be an element in $ N(X). $ Similarly, there is a moduli space $ M_H(X,P) $ (resp. $ M_H^s(X,x) $) of semistable (resp. stable) sheaves on $ X $ with numerical class $ x. $
\begin{rmk}
	The moduli space $ M_H(X, x) $ may be empty. If we choose $ x = -1 = -[\cO_X]$ in $ N(X), $ then the moduli space $ M_H(X, -1) $ is empty for an obvious reason. This perhaps is a hint that we should consider moduli spaces of semistable complexes of coherent sheaves such that $ M_H(X,-1) $ would then consist of the 2-term complex $ 0 \to \cO_X $ concentrated in degree $ -1 $ and $ 0. $
\end{rmk}
Now we recall the slope and its corresponding $ \mu $-stability of a coherent sheaf on a smooth projective scheme.
\begin{defn}
	Let $ (X, H) $ be a connected smooth polarized projective scheme in Notation \ref{notation_rk_deg}. We define a \tb{slope map}\footnote{This is neither additive nor multiplicative.}
	\begin{equation*}
	\mu: N(X) \to \QQ \cup \{\infty\}
	\end{equation*}
	by setting
	\begin{equation*}
	\mu(x) = \left \{
	\begin{array}{ll}
	\frac{\deg(x)}{\rk(x)} & \mbox{if $\rk(x) \neq 0$};\\
	\infty & \mbox{if $\rk(x) = 0$}.
	\end{array}
	\right.
	\end{equation*}
	If $ h = c_1(H) $ in $ R^1(X), $ then we will use the notation $ \mu_h $ instead of $ \mu $ when we want to emphasize the dependence of the slope map on $ h. $
\end{defn}
\begin{defn}
	Let $ (X, H) $ be a connected smooth polarized projective scheme. We say a coherent sheaf $ E $ on $ X $ is \tb{$ \mu $-(semi)stable} if $ E $ is a nonzero torsion-free sheaf and
	\begin{equation*}
	\mu(F) \ (\leq) \ \mu(E)
	\end{equation*}
	for all subsheaves $ F $ of $ E $ with $ 0 < \rk(F) < \rk(E). $
	We also refer to \tb{$ \mu $-stability} as \tb{slope stability}.
\end{defn}

The notions of stability and $ \mu $-stability are equivalent for torsion-free sheaves on a smooth polarized projective curve. In general, we have the following
\begin{lem}[{\cite[Lemma 1.2.13 and Lemma 1.2.14]{huybrechts2010geometry}}]\label{lem_slope_Gieseker}
	Let $ (X, H) $ be a connected smooth polarized projective scheme. Let $ E $ be a nonzero torsion-free sheaf $ E $ on $ X $. Then we have a chain of implications:
	\begin{equation*}
	E \ \text{is}\ \mu\text{-stable}\ \Rightarrow\ E\ \text{is stable}\ \Rightarrow\ E\ \text{is semistable}\ \Rightarrow\ E \ \text{is}\ \mu\text{-semistable}.
	\end{equation*}
	If $ \gcd(\rk(E),\deg(E)) = 1,$ then we also have the implication
	\begin{equation*}
	E \ \text{is}\ \mu\text{-semistable}\ \Rightarrow\ E \ \text{is}\ \mu\text{-stable}.
	\end{equation*}
\end{lem}
\begin{rmk}
	Lemma \ref{lem_slope_Gieseker} tells that if $ x $ is an element in $ N(X) $ with 
	$$ \rk(x) > 0 \quad \text{and} \quad \gcd(\rk(x),\deg(x)) = 1, $$
	then $ M(X,x) = M^s(X,x). $ In particular, if $ \rk(x) = 1, $ then $ M(X,x) = M^s(X,x). $
\end{rmk}

	In general, we expect the moduli space $ M(X,x) $ or $ M^s(X,x) $ to consist of singular connected components of various dimensions even if $ X $ is a smooth projective variety. The situation is much better when $ X $ has nice geometric properties, for example, when $ X $ is a $ K3 $ surface.
Consider a $ K3 $ surface $ X $. Since the Mukai vector map $ v: N(X) \to R(X) $ is an isomorphism of abelian groups, fixing an element $ x $ in $ N(X) $ is the same as fixing a vector $ v $ in $ R(X) $ for moduli spaces of stable sheaves on $ X $. Recall a few notions.
\begin{defn}
	A vector $ v $ in $ R(X) $ is \tb{primitive} if it is not a multiple of another vector $ w $ in $ R(X) $. An element $ x $ in $ N(X) $ is said to be primitive if it is not a multiple of another element $ y $ in $ N(X). $
\end{defn}
\begin{defn}
	The \tb{ample cone} of $ X $ is a convex cone defined by
	$$ \Amp(X) = \{C \in R^1(X)_\RR: C = \sum a_i [D_i] \ \text{where } a_i > 0 \ \text{and} \ D_i \ \text{is an ample divisor} \}. $$
\end{defn}
It's shown in \cite{yoshioka1994chamber} that a vector $ v $ in $ R(X) $ determines a countable locally finite set of hyperplanes in $ \Amp(X) $, which are called \tb{$ v $-walls}.
\begin{defn}
	Let $ x $ be an element in $ N(X) $ with Mukai vector $ v $ in $ R(X). $ An ample line bundle $ H $ on $ X $ is said to be \tb{$ x $-generic}, or \tb{$ v $-generic}, if the divisor class $ h = c_1(H) \in \Amp(X) $ does not lie on any of the $ v $-walls.
\end{defn}
Under some conditions, $ \mu $-semistability implies $ \mu $-stability.
\begin{lem}\label{lem_slope_semistable_implies_stable}
	Let $ (X, H) $ be a polarized $ K3 $ surface. Let $ x $ be an element in $ N(X) $ with Mukai vector $ v = (r, c_1, s) $ in $ R(X) $ with $ r > 0 $. Suppose $ (r, c_1) $ is primitive and $ H $ is $ v $-generic. If a torsion-free sheaf $ E $ on $ X $ with numerical class $ x $ is $ \mu $-semistable, then it is $ \mu $-stable. 
\end{lem}
\begin{proof}
	Take a $ \mu $-semistable sheaf $ E $ on $ X $ with numerical class $ x. $ Suppose it is not $ \mu $-stable. Then there is a subsheaf $ F \subset E $ with $ 0 < \rk(F) < \rk(E) $ such that $ \mu(F) = \mu(E), $ i.e.,
	\begin{equation*}
		\frac{\deg(F)}{\rk(F)} = \frac{\deg(E)}{\rk(E)}.
	\end{equation*}
	Since $ H $ is $ x $-generic, we must have $ \rk(E) c_1(F) = \rk(F) c_1(E) $ in $ R^1(X). $ This implies $ (r, c_1) = r/\rk(F)(\rk(F),c_1(F)), $ which is impossible since $ (r, c_1) $ is primitive.
\end{proof}
In general, fixing a primitive element $ x $ in $ N(X) $ such that the polarization is $ x $-generic guarantees that semistability is the same as stability. We combine the results in \cite[Corollary 4.6.7]{huybrechts2010geometry}, \cite[Corollary 10.2.1]{huybrechts2016lectures}, and \cite[Proposition 10.2.5]{huybrechts2016lectures} into the following
\begin{prop}\label{prop_K3_smoothness_and_dim}
	Let $ (X, H) $ be a polarized $ K3 $ surface. Let $ x $ be an element in $ N(X) $ with Mukai vector $ v  = (r,c_1,s) $ in $ R(X) $ and degree $ d = \deg(x) $. Then $ M_H^s(X,x) $ is either empty or a smooth quasi-projective scheme of dimension $ n = 2 - \inprod{v}^2. $ Suppose the Mukai vector $ v  $ satisfies either of the following conditions:
	\begin{enumerate}[font=\normalfont,leftmargin=2em]
		\item $ \gcd(r, d, s) = 1,$
		\item $ x $ is primitive and $ H $ is $ x $-generic.
	\end{enumerate}
	Then $ M_H(X, x) = M_H^s(X, x), $ and it is either empty or a smooth projective scheme of dimension $ n. $ Moreover, $ M_H(X, x) $ is a fine moduli space in case $ (1). $
\end{prop}
Recall that the Hilbert scheme $ \Hilb^n(X) $ for a $ K3 $ surface $ X $ is an irreducible symplectic manifold of dimension $ 2n, $ and $ \Hilb^n(X) $ is deformation equivalent to $ \Hilb^n(Y) $ for any other $ K3 $ surface $ Y. $ Now we can state the following classical result.
\begin{thm}[{\cite[Theorem 8.1]{yoshioka2001moduli}}]\label{thm_yoshioka}
	Let $ (X, H) $ be a polarized $ K3 $ surface. Let $ x $ be an element in $ N(X) $ with $ \rk(x) > 0 $. Assume $ x $ is primitive and $ H $ is $ x $-generic. Then $ M_H(X, x) = M_H^s(X, x), $ which is non-empty if and only if $ \inprod{v(x)^2} \leq 2. $
	If $ M_H(X, x) $ is non-empty, then it is an irreducible symplectic manifold of dimension $ n = 2 - \inprod{v^2} $ deformation equivalent to $ \Hilb^{n/2}(X) $.
\end{thm}
We will generalize this result to equivariant moduli spaces of stable sheaves on a $ K3 $ surface with symplectic automorphisms in Section \ref{sec_proof_main_thm}.
\vspace{4pt}

\subsection{Moduli spaces of sheaves on projective stacks}\label{sec_moduli_sheaves_stacks2}
In this section we review Gieseker stabilities of coherent sheaves on projective stacks over an algebraically closed field in \cite{nironi2009moduli}. 


Let $ S $ be a scheme. We first recall the notion of tame stacks.
\begin{defn-prop}[{\cite[Definition 3.1 and Theorem 3.2]{abramovich2008tame}}]
	Let $ \cX $ be an Artin stack locally of finite presentation over $ S $ with a finite inertia, i.e., the canonical morphism $ I\cX \to \cX $ is finite. Let $ \pi: \cX \to Y $ denote the coarse moduli space of $ \cX $. The stack $ \cX $ is \tb{tame} if it satisfies either of the following equivalent conditions:
	\begin{enumerate}[font=\normalfont,leftmargin=*]
		\item The pushforward functor $ \pi_*: \Qcoh(\cX) \to \Qcoh(Y) $ is exact.
		\item If $ k $ is an algebraically closed field with a morphism $ \Spec k \to S $ and $ \eta $ is an object of $ \cX(k) $, then the automorphism group $ \Aut_k(\eta) $ is linearly reductive over $ k $, i.e., the order of $ \Aut_k(\eta) $ is coprime to $ \characteristic k. $
	\end{enumerate}
\end{defn-prop}

Let $ S = \Spec k $ for a field $ k. $ Recall the notion of generating sheaves on tame DM stacks over $ k $ introduced in \cite{olsson2003quot}. 
\begin{defn-prop}[{\cite[Theorem 5.2]{olsson2003quot}}]
	Let $ \cX $ be a tame DM stack over $ k $ with a coarse moduli space $ \pi: \cX \to Y $. A locally free sheaf $ \cV $ on the stack $ \cX $ is a \tb{generating sheaf} if it satisfies either of the following equivalent conditions:
	\begin{enumerate}[font=\normalfont,leftmargin=*]
		\item For every geometric point $ x $ in $ \cX, $ the local representation $ \phi_x: G_x \to \GL(V_x) $ of the automorphism group $ G_x $ on the fiber $ V_x $ of $ \cV $ at $ x $ contains every irreducible representation of $ G_x $.
		\item For every quasi-coherent sheaf $ \cF $ on $ \cX, $ the natural morphism
		\begin{equation*}
		\left(\pi^*\pi_*(\cF \otimes \cV^\vee)\right) \otimes \cV \to \cF
		\end{equation*}
		is surjective.
	\end{enumerate}
\end{defn-prop}
\begin{ex}
	Let $ \cX = BG $ for a finite group $ G $. Then the regular representation $ \rho_\reg $ of $ G $ is a generating sheaf on $ BG. $
\end{ex}
\begin{ex}
	For a quotient stack $ \cX = [X/G] $ where $ G $ is a finite group, the vector bundle $ \cO_\cX \otimes \rho_\reg $ is a generating sheaf on $ \cX $. We will consider this example again in the next section.
\end{ex}
\begin{rmk}
	In general we do not know whether there are generating sheaves on tame DM stacks. However, they do exist on projective stacks. Generating sheaves are a key ingredient in the construction of moduli spaces of sheaves on projective stacks.
\end{rmk}

\begin{defn}[{\cite[Definition 2.20]{nironi2009moduli}}]\label{defn_proj_stack}
	Let $ \cX $ be a separated quotient stack over $ k $. We say $ \cX $ is a \tb{projective} (resp. quasi-projective) stack if it is a tame DM stack and its coarse moduli space is a projective (resp. quasi-projective) scheme.
\end{defn}
Now we assume the base field $ k $ is algebraically closed. Let $ \cX $ be a projective stack over $ k $ with a coarse moduli space $ \pi: \cX \to Y. $ Since $ \pi: \cX \to Y $ is a proper morphism and $ \cX $ is tame, we have an exact functor
\begin{equation*}
\pi_*: \Coh(\cX) \to \Coh(Y).
\end{equation*}
Therefore, for every coherent sheaf $ \cE $ on $ \cX $, its pushforward $ \pi_*\cE $ is a coherent sheaf on $ Y $ and we have
\begin{equation*}
H^i(\cX, \cE) = H^i(Y, \pi_*\cE)
\end{equation*}
for all $ i \geq 0 $. Furthermore, since $ \cX $ is a quasi-compact quotient stack, it has a generating sheaf by \cite[Theorem 5.5]{olsson2003quot}.  Since $ Y $ is a projective scheme, it has ample line bundles. To facilitate further discussions, we make the following
\begin{defn}
	Let $ \cX $ be a projective stack over $ k $ with a coarse moduli space $ \pi: \cX \to Y. $ 
	Let $ \cV $ be a generating sheaf on $ \cX $, and let $ L $ be an ample line bundle on $ Y. $ Let $ \cH = \pi^*L. $ The pair $ (\cH, \cV) $ is called a \tb{polarization} of $ \cX. $ The triple $ (\cX, \cH, \cV) $ is called a \tb{polarized projective stack}.
\end{defn}
Let $ \cE $ be a nonzero coherent sheaf on $ \cX. $ The dimension of $ \cE $ is defined to be the dimension of its support, as in the case of schemes. By \cite[Proposition 3.6]{nironi2009moduli}, we have
\begin{equation*}
\dim \cE = \dim \pi_*(\cE \otimes \cV^\vee).
\end{equation*}
By the projection formula for DM stacks, we have
\begin{equation*}
\chi(\cX, \cE \otimes \cV^\vee \otimes \cH^{\otimes m}) = \chi(Y, \pi_*(\cE \otimes \cV^\vee) \otimes L^{\otimes m}) = P_{\pi_*(\cE \otimes \cV^\vee)}(m)
\end{equation*}
in $ \ZZ $ for all $ m $ in $ \ZZ. $ Therefore, the assignment $ m \mapsto \chi(\cX, \cE \otimes \cV^\vee \otimes \cH^{\otimes m}) $ is a numerical polynomial of degree $ n = \dim \cE. $ Now we can define a modified Hilbert polynomial of the sheaf $ \cE $ on the stack $ \cX $ in a similar fashion to the usual Hilbert polynomial of a sheaf on schemes.
\begin{defn-prop}
	Let $ (\cX, \cH, \cV) $ be a polarized projective stack over $ k $. Let $ \cE $ be a coherent sheaf on $ \cX. $ The \tb{modified Hilbert polynomial} of $ \cE $ with respect to the polarization $ (\cH, \cV) $ is the unique numerical polynomial
	\begin{equation*}
	\orbP_\cE(z) = a_n(\cE) z^n + \cdots + a_1(\cE) + a_0(\cE)
	\end{equation*}
	in $ \QQ[z] $ such that
	\begin{equation*}
	\orbP_\cE(m) = \chi(\cX, \cE \otimes \cV^\vee \otimes \cH^{\otimes m})
	\end{equation*}
	in $ \ZZ $ for all $ m $ in $ \ZZ. $ If $ \cE = 0, $ then $ \orbP_\cE(z) = 0; $ else, $ n = \dim \cE. $
	If $ \cE \neq 0, $ then the reduced Hilbert polynomial of $ \cE $ is defined by the monic polynomial
	\begin{equation*}
	\orbp_\cE(z) = \frac{\orbP_\cE(z)}{a_n(\cE)}
	\end{equation*}
	in $ \QQ[z]. $
\end{defn-prop}
\begin{rmk}
	The modified Hilbert polynomial is additive on short exact sequences on $ \Coh(\cX) $ because the functor $ \Coh(\cX) \to \Coh(Y), \ \cE \to \pi_*(\cE \otimes \cV^\vee) $ is a composition of exact functors and hence is exact. Therefore, the modified Hilbert polynomial descends to an additive map
	\begin{equation*}
	\orbP: N(\cX) \to \QQ[z].
	\end{equation*}
\end{rmk}
We define pure sheaves on projective stacks as usual.
\begin{defn}
	A nonzero coherent sheaf $ \cE $ on $ \cX $ is pure if $ \dim \cF = \dim \cE $ for all nonzero subsheaves $ \cF \subset \cE $ on $ \cX. $ 
\end{defn}
\begin{defn}
	Let $ (\cX, \cH, \cV) $ be a polarized projective stack over $ k $. A coherent sheaf $ \cE $ on $ \cX $ is \tb{$ (\cH, \cV) $-(semi)stable} if it is a nonzero pure sheaf and
	$$ \orbp_\cF(z) \ \ (\leq)\ \ \orbp_\cE(z) $$ 
	for all proper nonzero subsheaves $ \cF \subset \cE $. When the choice of $ \cV $ is understood, $ (\cH, \cV) $-stability is also called \tb{$ \cH $-stability}; if both $ \cH $ and $ \cV $ are understood, $ (\cH, \cV) $-stability is also called \tb{stability}.
\end{defn}
Fix a numerical polynomial $ P $ in $ \QQ[z] $. By a GIT construction of Nironi in \cite{nironi2009moduli}, there is a projective scheme $ M_{\cH,\cV}(\cX,P) $ parametrizing $ S $-equivalent classes of $ (\cH,\cV) $-semistable sheaves on $ \cX $ with modified Hilbert polynomial $ P, $ which compactifies a quasi-projective scheme $ M_{\cH,\cV}^s(\cX,P) $ parametrizing $ (\cH,\cV) $-stable sheaves on $ \cX $ with modified Hilbert polynomial $ P. $ Similarly, choosing a numerical class $ x $ in the numerical $ K $-theory $ N(\cX) $ gives a moduli space $ M_{\cH,\cV}(\cX,x) $ (resp. $ M_{\cH,\cV}^s(\cX,x) $) of $ (\cH,\cV) $-semistable (resp. $ (\cH,\cV) $-stable) sheaves on $ \cX $ with numerical class $ x. $

Now we generalize a smoothness result of moduli spaces of sheaves on projective schemes in \cite[Proposition 10.1.11]{huybrechts2016lectures} to projective stacks.
\begin{prop}\label{prop_smoothness_criterion}
	Let $ (\cX, \cH, \cV) $ be a smooth polarized projective stack over $ k $. Choose an element $ x $ in $ N(\cX) $ such that the moduli space $ M = M_{\cH,\cV}(\cX,x) $ is non-empty. Let $ t $ be a closed point in $ M $ corresponding to a stable sheaf $ \cE $ in the groupoid $ \cM(k) $.
	\begin{enumerate}[font=\normalfont,leftmargin=*]
		\item There is a natural isomorphism $ T_{t} M \cong \Ext^1(\cE,\cE) $.
		\item If $ \Ext^2(\cE,\cE) = 0, $ then $ M $ is smooth at $ t. $
		\item If $ \Pic(\cX) $ is smooth and the trace map $ \Ext^2(\cE,\cE) \to H^2(\cX, \cO_\cX) $ is injective, i.e., $ \Ext^2(\cE,\cE)_0 = 0, $ then $ M $ is smooth at $ t. $
	\end{enumerate}
\end{prop}
\begin{proof}
	Let $ M^s = M_{\cH,\cV}^s(\cX, x). $ The element $ x $ in $ N(\cX) $ determines an integer $ r $ such that all semistable sheaves with modified Hilbert polynomial $ \orbP_x $ are $ r $-regular. In particular, $ \cE $ is $ r $-regular, which means $ \pi_*(\cE \otimes \cV^\vee) $ is $ r $-regular. Let $ N = \orbP_x, $ and let $ \cF = \cV^{\oplus {N}} \otimes \cH^{\otimes -r}.$
	By Corollary 4.71 in \cite{chen2023orbifold}, there is an open subscheme $ R^s \subset \Quot_\cX(\cF,\orbP_x) $ such that the geometric quotient
	\begin{equation*}
	\pi: R^s \to M^s
	\end{equation*}
	is a principal $ \PGL(N,k) $-bundle. Consider a point $ q $ in the fiber $ \pi^{-1}(t), $ which corresponds to a short exact sequence
	\begin{equation}\label{eq_quot_ses}
	0 \to \cK \xrightarrow{s} \cF \xrightarrow{q} \cE \to 0
	\end{equation}
	in $ \Coh(\cX) $. Let $ \cO_q = \PGL(N,k) \cdot q $ denote the orbit of $ q $ in $ R^s. $ Then we have
	\begin{equation*}
	T_t M \cong T_q R^s/T_q \cO_q,
	\end{equation*}
	where the tangent space $ T_q R^s \cong \Hom(\cK, \cE) $. Now we apply $ \Hom(\ {\cdot} \ ,\cE) $ to (\ref{eq_quot_ses}) and get a long exact sequence
	\begin{equation*}
	0 \to \End(\cE) \to \Hom(\cF,\cE) \xrightarrow{s^\sharp} \Hom(\cK,\cE) \to \Ext^1(\cE,\cE) \to \Ext^1(\cF,\cE) \to \cdots
	\end{equation*}
	We claim that 
	$$ \Ext^i(\cF,\cE) = 0 \quad \text{for all}\ i \geq 1. $$
	This is a consequence of the $ r $-regularity of $ \cE. $ Take an integer $ i \geq 1. $ Since $ \cF $ is locally free, we have:
	\begin{align*}
	\Ext^i(\cF,\cE) & = H^i(\cX, \cE \otimes \cF^\vee) \\
	& = H^i(\cX, (\cE \otimes \cV^\vee \otimes \cH^r)^{\oplus N}) \\
	& = \bigoplus_{k=1}^N H^i(\cX, \cE \otimes \cV^\vee \otimes \cH^{\otimes r}) \\
	& = \bigoplus_{k=1}^N H^i(Y, \pi_*(\cE \otimes \cV^\vee) \otimes L^{\otimes r}) & \text{since $ \cX $ is tame}\\
	& = 0 & \text{since $ \pi_*(\cE \otimes \cV^\vee) $ is $ r $-regular}\\
	\end{align*}
	The claim is proved. Note that $ T_q \cO_q \cong \img(s^\sharp) \subset \Hom(\cK, \cE). $  Therefore, we have
	\begin{equation*}
	\Ext^1(\cE,\cE) \cong \Hom(\cK, \cE)/ \img(s^\sharp) \cong T_t M.
	\end{equation*}
	From the vanishing of $ \Ext^2(\cF, \cE), $ we also have
	\begin{equation*}
	\Ext^1(\cK,\cE) \cong \Ext^2(\cE, \cE),
	\end{equation*}
	where $ \Ext^1(\cK,\cE) $ is the space of obstructions to deform $ \cE $ in $ M. $ Now it remains to prove the last statement. Since $ \cX $ is smooth and has the resolution property, there is a determinant morphism
	\begin{equation*}
	\det: M \to \Pic(\cX), \quad \cE \to \det(\cE),
	\end{equation*}
	where $ \det(\cE) $ is well-defined via a finite locally resolution of $ \cE. $ Since $ \det(\cE) $ is a line bundle, we have an isomorphism of $ k $-vector spaces:
	\begin{equation*}
	H^2(\cX, \cO_\cX) \cong \Ext^2(\det(\cE),\det(\cE)).
	\end{equation*}
	Therefore, we can identify the trace map as the following map
	\begin{equation*}
	\tr: \Ext^2(\cE,\cE) \to \Ext^2(\det(\cE),\det(\cE)).
	\end{equation*}
	Consider a local Artinian $ k $-algebra $ (A, \frak{m}) $ where $ \frak{m} $ is the maximal ideal in $ A $. Let $ I $ be an ideal in $ A $ such that $ I \cdot \frak{m} = 0.$ Then the trace map induces a map
	\begin{equation*}
	\tr \otimes_k I: \Ext^2(\cE,\cE) \otimes_k I \to \Ext^2(\det(\cE),\det(\cE)) \otimes_k I
	\end{equation*}
	which sends the obstruction $ \ob(\cE,A) $ to lift an $ A/I $-flat deformation of $ \cE $ in $ M $ to the obstruction $ \ob(\det(\cE), A) $ to lift an $ A/I $-flat deformation of $ \det(\cE) $ in $ \Pic(\cX) $. Suppose the trace map is injective. Then the $ k $-linear map $ \tr \otimes_k I $ is injective. Therefore, the smoothness of $ \Pic(\cX) $ implies $ \ob(\det(\cE), A) = 0 $, and hence $ \ob(\cE,A) = 0. $
\end{proof}

\subsection{Equivariant moduli spaces of sheaves on projective schemes}
In this section we consider the following
\begin{situation}\label{situation_finite_group}
	Let $ X $ be a projective scheme over $ \CC $, and let $ G $ be a finite group acting on $ X $. Then we have a commutative triangle
	\vspace{5pt}
	\begin{equation*}
	\begin{tikzcd}[column sep=1.6em,row sep=3.2em]
	& X \arrow{ld}[swap]{p} \arrow{rd}{q} & \\
	\cX \arrow{rr}{\pi} & & Y,
	\end{tikzcd}\vspace{10pt}
	\end{equation*}
	where $ p $ and $ q $ are finite, $ \pi $ is proper, $ \cX = [X/G] $ is a projective stack, $ Y = X/G $ is a projective scheme which is both the geometric quotient of $ X $ by $ G $ and the coarse moduli space of $ \cX $. Let $ L $ be an ample line bundle on $ Y $.
	Put
	$$ \cH = \pi^*L \quad \text{and} \quad H = q^*L. $$
	Then $ H $ is an ample line bundle on $ X $ since it's a pullback of an ample line bundle along a finite morphism.
	The line bundle $ \cH $ on $ \cX $ corresponds to an $ G $-equivariant ample line bundle $ (H, \phi) $ on $ X. $ Since it's a pullback from a line bundle on the coarse moduli space $ M, $ all local representations 
	\begin{equation*}
	\phi_x: G_x \to \GL(H_x)
	\end{equation*}
	of the stabilizers $ G_x $ on the fibers of $ H $ at every point $ x \in X $ are trivial.
	Consider the regular representation $ \rho_\reg $ of $ G $. Pulling it back along $ \cX \to BG $ gives a canonical generating sheaf 
	$$ \cV_\reg = \cO_\cX \otimes \rho_\reg $$
	on the stack $ \cX. $ Then we have a polarized projective scheme $ (X, H) $ and a polarized projective stack $ (\cX, \cH, \cV_{\reg}) $.
\end{situation}
By the categorical equivalence $ \Coh(\cX) \cong \Coh^G(X), $ we define the stability of $ G $-equivariant sheaves on $ X $ as the stability of the corresponding sheaves on $ \cX. $
\begin{defn}\label{stability_equivariant_sheaves}
	Let $ \cX = [X/G] $ in Situation \ref{situation_finite_group}. A $ G $-equivariant sheaf $ (E,\phi) $ on $ X $ is \tb{$ \cH $-semistable }(resp. \tb{$ \cH $-stable}) if the corresponding sheaf $ \cE $ on $ \cX $ is $ (\cH, \cV_\reg) $-semistable (resp. $ (\cH, \cV_\reg) $-stable).
\end{defn}
Fix an element $ x $ in $ N(\cX). $ There is a projective scheme $ M_\cH(\cX,x) $ parametrizing $ S $-equivalent classes of Gieseker $ (\cH, \cV_\reg) $-semistable sheaves on $ \cX $, which compactifies the quasi-projective scheme $ M_\cH^s(\cX,x) $ parametrizing Gieseker $ (\cH, \cV_\reg) $-stable sheaves on $ \cX $. These will be $ G $-equivariant moduli spaces of sheaves on $ X. $
\begin{defn}
	Let $ \cX = [X/G] $ in Situation \ref{situation_finite_group}. Let $ x $ be an element in $ N(\cX). $ We call the projective (resp. quasi-projective) scheme $ M_\cH(\cX,x) $ (resp. $ M_\cH^s(\cX,x) $) the \tb{$ G $-equivariant moduli space} of $ \cH $-semistable (resp. $ \cH $-stable) sheaves on $ X $ with numerical class $ x $.
\end{defn}

Let $ \cX = [X/G] $ in Situation \ref{situation_finite_group}. The sheaf cohomologies on $ X $ and $ \cX $ are related as follows.
\begin{lem}
	For a coherent sheaf $ \cE = (E, \phi) $ on $ \cX $ and an integer $ m, $
	\begin{equation*}
	H^i(X, E \otimes H^{\otimes m}) \cong H^i(\cX, \cE \otimes \cV_\reg \otimes \cH^{\otimes m})
	\end{equation*}
	for all integer $ i \geq 0, $ and hence
	\begin{equation*}
	\chi(X, E \otimes H^{\otimes m}) = \chi(\cX, \cE \otimes \cV_\reg \otimes \cH^{\otimes m}).
	\end{equation*}
\end{lem}
\begin{proof}
	Take a coherent sheaf $ \cE = (E, \phi) $ on $ \cX. $ Fix an integer $ m $ and an integer $ i \geq 0. $ We first show\footnote{The argument here is due to Promit Kundu.} $ H^i(X, E \otimes H^{\otimes m}) \cong H^i(\cX, \cE \otimes \cV_\reg \otimes \cH^{\otimes m}). $  
	Since $ p: X \to \cX $ is finite and affine, the pushforward functor $ p_*: \Coh(X) \to \Coh(\cX) $
	is exact. Note that $ p_*\cO_X \cong \cO_\cX \otimes \rho_\reg = \cV_\reg $.  Applying the projection formula to $ p: X \to \cX $, we have
	\begin{equation*}
	p_*p^*(\cE \otimes \cH^{\otimes m}) \cong \cE \otimes \cH^{\otimes m} \otimes p_*\cO_X \cong \cE \otimes \cH^{\otimes m} \otimes \cV_\reg,
	\end{equation*}
	and hence
	\begin{align*}
	H^i(X, E \otimes H^{\otimes m}) \cong H^i(X, p^*(\cE \otimes \cH^{\otimes m})) \cong H^i(\cX, \cE \otimes \cV_\reg \otimes \cH^{\otimes m}).
	\end{align*}
\end{proof}
This lemma identifies the modified Hilbert polynomials for sheaves $ \cE $ on $ \cX $ with the Hilbert polynomials for the pullback $ E = p^*\cE $ on $ X $.
\begin{cor}\label{cor_P_orb_equals_P}
	Let $ \cE = (E, \phi) $ be a coherent sheaf on $ \cX $. Then
	$$ \orbP_\cE(z) = P_E(z) $$
	in $ \QQ[z] $.
\end{cor}
\begin{proof}
	Note that $ \rho_\reg $ is self-dual and hence $ \cV_\reg \cong \cV_\reg^\vee. $ We then have
	\begin{align*}
	P_E(m) = \chi(X, E \otimes H^{\otimes m}) = \chi(\cX, \cE \otimes \cV_\reg \otimes \cH^{\otimes m}) = \orbP_\cE(m)
	\end{align*}
	in $ \ZZ $ for all $ m $ in $ \ZZ $.
\end{proof}
\begin{rmk}
	An alternative definition of stability of $ G $-equivariant sheaves on $ X $ was introduced in \cite{amrutiya2015moduli} as follows: a $ G $-equivariant sheaf $ (E, \phi) $ on $ X $ is (semi)stable if $ E $ is a nonzero pure sheaf and $ p_F(z) \ (\leq)\ p_E(z) $ for all nonzero proper $ G $-equivariant subsheaves $ (F, \phi) \subset (E, \phi) $. Corollary \ref{cor_P_orb_equals_P} implies that this definition is equivalent to Definition \ref{stability_equivariant_sheaves}.
\end{rmk}

The canonical morphism $ p: X \to \cX $ is smooth and hence induces an exact pullback functor 
$$ p^*: \Coh(\cX) \to \Coh(X), \quad \cE = (E, \phi) \mapsto E, $$ 
which has a left adjoint 
$$ \Ind: \Coh(X) \to \Coh(\cX) $$ 
defined as follows. Let $ E $ be a coherent sheaf on $ X. $ It determines a $ G $-equivariant coherent sheaf
\begin{equation*}
\Ind E = \left(S_G(E), \phi \right), \ \text{with}\ S_G(E) = \bigoplus_{g \in G} \{g\} \times g^*E.
\end{equation*}
Here $ \phi $ is a canonical $ G $-equivariant structure on $ S_G(E) $ such that for each $ x \in G $, and each $ h \in G, $ there is a linear isomorphism
\begin{align*}
\phi_{x,h}: \bigoplus_{g \in G} \{g\} \times E_{gx} & \to \bigoplus_{g \in G} \{g\} \times E_{ghx} \\
(g,v_g)_{g \in G} & \mapsto (gh^{-1},v_g)_{g \in G} = (g,v_{gh})_{g \in G}.
\end{align*}
Note that $ \phi_{x,h} $ is a permutation if $ E $ is locally free.
We say the sheaf $ \Ind E $ is induced from $ E $. 
Now we have the following
\begin{lem}\label{lem_induced_sheaves}
	Let $ E $ be a coherent sheaf on $ X $ which induces a sheaf $ \cE = \Ind E = (S_G(E), \phi) $ on $ \cX $. Let $ \cF = (F, \psi) $ be a coherent sheaf on $ \cX. $
	Then there is an isomorphism of vector spaces
	\begin{equation*}
	\Hom_{\Coh(\cX)}(\cE, \cF) \xrightarrow{\sim} \Hom_{\Coh(X)}(E, F).
	\end{equation*}
\end{lem}
\begin{proof}
	The sheaf $ \cF $ corresponds to a $ G $-equivariant sheaf $ \cF = (F, \psi) $ on $ X. $ 
	A morphism $ \lambda: \Ind E \to \cF $ pulls back to a morphism $ p^*\lambda: S_G(E) \to F, $ which induces a morphism $ f: E \to F $ given by the composition
	\begin{equation*}
	f: E \into S_G(E) \xrightarrow{p^*\lambda} F.
	\end{equation*}
	On the other hand, suppose we have a morphism $ f: E \to F. $ Consider any point $ x $ in $ X $ and an element $ g $ in $ G. $ There are fiberwise linear isomorphisms
	\begin{equation*}
	f_{gx}: E_{gx} \xrightarrow{\sim} F_{gx} \quad \text{and} \quad \psi_{x,g}: F_x \xrightarrow{\sim} F_{gx}.
	\end{equation*}
	Define a morphism $ \lambda: S_G(E) \to F $ fiberwise by
	\begin{align*}
	\lambda_x: \bigoplus_{g \in G} \{g\} \times E_{gx} & \to F_x \\
	(g,v_g)_{g \in G} & \mapsto \sum_{g \in G} \psi_{x,g}^{-1} \left(f_{gx}(v_g)\right). 
	\end{align*}
	Take any element $ h $ in $ G. $ Then a direct computation show that
	\begin{equation*}
	\lambda_{hx} \circ \phi_{x,h} = \psi_{x,h} \circ \lambda_x,
	\end{equation*}
	which means $ \lambda $ is a $ G $-equivariant structure on $ S_G(E). $ By construction, the assignment $ \lambda \mapsto f $ and $ f \mapsto \lambda $ are inverses to each other.
\end{proof}
An immediate consequence of Lemma \ref{lem_induced_sheaves} is the following
\begin{prop}\label{prop_Gieseker_relation}
	Let $ \cX = [X/G] $ in Situation \ref{situation_finite_group}. Let $ \cE = (E, \phi) $ be a nonzero coherent sheaf on $ \cX $. Then $ \cE $ is $ \cH $-semistable if and only if $ E $ is $ H $-semistable. If $ E $ is $ H $-stable, then $ \cE $ is $ \cH $-stable.
\end{prop}
\begin{proof}
	By Remark 3.3 in \cite{nironi2009moduli}, $ \cE $ is pure of dimension $ d $ if and only if $ E $ is so. Assume $ \cE $ is a pure sheaf, otherwise there is nothing to prove. Suppose $ \cE $ is not (semi)stable. Then there is a proper nonzero subsheaf $ \cF \subset \cE $ such that
	\begin{equation*}
	\orbp_\cF(z)\ \ (\geq) \ \ \orbp_\cE(z)
	\end{equation*}
	in $ \QQ[z] $. Since $ p: X \to \cX $ is a flat morphism, the pullback functor $ p^*: \Coh(\cX) \to \Coh(X) $ is exact, which implies $ F = p^*\cF \subset E $ is also a proper nonzero subsheaf. By Corollary \ref{cor_P_orb_equals_P}, $ p_F(z) (\geq) p_E(z), $ therefore $ E $ is not (semi)stable. Hence we have proved two implications:
	\begin{equation*}
	\text{$ E $ is (semi)stable} \quad \Rightarrow \quad \text{$ \cE $ is (semi)stable}.
	\end{equation*}
	Suppose $ E $ is not semistable. Then there is a proper nonzero subsheaf $ F \subset E $ such that
	\begin{equation*}
	p_F(z) > p_E(z)
	\end{equation*}
	in $ \QQ[z]. $ Moreover, by \cite[Lemma 1.3.5]{huybrechts2010geometry}, $ F $ can be chosen as a semistable sheaf. Let $ \cF = \Ind F = (S_G(F), \phi). $ By Lemma \ref{lem_induced_sheaves}, the inclusion $ F \into E $ gives a morphism 
	$$ \lambda: \cF \to \cE $$
	whose image is a subsheaf $ \cE' \subset \cE $ which is nonzero since $ F \into E $ is not a zero morphism. Since $ \cF $ is semistable and $ \cF \to \cE' $ is a quotient, we have
	\begin{equation*}
	\orbp_{\cE'}(z) \ \geq\ \orbp_\cF(z) = p_{S_G(F)}(z)
	\end{equation*}
	in $ \QQ[z]. $ Take an element $ g \in G. $ The morphism
	\begin{equation*}
	g: X \to X, \quad x \mapsto gx
	\end{equation*}
	is smooth and is in particular flat, so $ g^*: \Coh(X) \to \Coh(X) $ is an exact functor. Therefore, we have
	\begin{equation*}
	H^0(X, F) = \Hom(\cO_X, F) \cong \Hom(\cO_X, g^*F) = H^0(X, g^*F),
	\end{equation*}
	which implies
	\begin{equation*}
	P_F(z) = P_{g^*F}(z) \quad \text{for all } g \in G.
	\end{equation*}
	By the definition of $ S_G(F), $ we then have
	\begin{equation*}
	P_{S_G(F)}(z) = n P_F(z),
	\end{equation*}
	where $ n = |G|, $ and hence 
	$$ p_{S_G(F)}(z) = p_F(z). $$ 
	Now we conclude
	\begin{equation*}
	\orbp_{\cE'}(z) \geq p_F(z) > p_E(z) = \orbp_\cE(z),
	\end{equation*}
	which shows $ \cE $ is not semistable. In other words, we have proved the implication
	\begin{equation*}
	\text{$ \cE $ is semistable} \quad \Rightarrow \quad \text{$ E $ is semistable}.
	\end{equation*}
\end{proof}

\begin{rmk}
	\begin{enumerate}[font=\normalfont,leftmargin=*]
		\item Let's summarize the results in Proposition \ref{prop_Gieseker_relation}:
		\vspace{5pt}
		\begin{equation}\label{diag_Gieseker_relation}
		\begin{tikzcd}[column sep=2.4em,row sep=3.6em]
		\text{$ E $ is stable} \arrow[r, Rightarrow, start anchor = {[xshift = 3.3ex]}, end anchor = {[xshift = -3.3ex]}] \arrow[d, Rightarrow, start anchor = {[yshift = -2.5ex]}, end anchor = {[yshift = 2.5ex]}] & \text{$ \cE $ is stable} \arrow[d, Rightarrow, start anchor = {[yshift = -2.5ex]}, end anchor = {[yshift = 2.5ex]}] \\
		\text{$ E $ is semistable} \arrow[r, Leftrightarrow, start anchor = {[xshift = 1ex]}, end anchor = {[xshift = -1ex]}] & \text{$ \cE $ is semistable}
		\end{tikzcd}\vspace{10pt}
		\end{equation}
		In general, the stability of a sheaf $ \cE = (E, \phi) $ on $ \cX $ is weaker than the stability of the sheaf $ E $ on $ X $. For example, consider $ G = \{1,\sigma\} \cong \mu_2 $ acting on $ X $ and let $ \cE = \cO_p $ be the structure sheaf of a generic point $ p \in \cX $. So $ E = p^*\cE = \cO_{Gx} $ where $ Gx = \{x, \sigma \cdot x\} $ is a free orbit on $ X $, i.e., $ E = \cO_x \oplus \cO_{\sigma \cdot x} $. Note that the $ G $-equivariant structure $ \phi $ on $ E $ is from the regular representation of $ G $. Therefore, $ \cE $ is stable on $ \cX $, but $ E $ is strictly semistable on $ X $.
		\item Let $ [X/G] $ be a connected smooth quotient stack in Situation \ref{situation_finite_group}. Suppose $ (E, \phi) $ is a $ G $-equivariant sheaf on $ X $ with $ \gcd(\rk(E),\deg(E)) = 1. $ By Lemma \ref{lem_slope_Gieseker}, the stability of $ E $ is the same as the semistability of $ E. $ In this case, every arrow in diagram (\ref{diag_Gieseker_relation}) becomes a two-sided arrow.
	\end{enumerate}
\end{rmk}
\vspace{5pt}

\section{Equivariant moduli theory on $ K3 $ surfaces}\label{sec_moduli_sheaves_K3/G}
\vspace{5pt}
In this section we will study equivariant moduli spaces of sheaves on $ K3 $ surfaces with symplectic automorphisms. We first introduce a notation.
\begin{notation}
	Let $ X $ be a $ K3 $ surface. Let $ G $ be a finite subgroup of the symplectic automorphisms of $ X $. Then the quotient stack $ \cX = [X/G] $ is a connected smooth projective stack over $ \CC. $ The quotient stack $ \cX $ is called a model of $ [K3/G] $ and is denoted by $ \cX = [K3/G] $.
\end{notation}
Take a model $ \cX = [K3/G] $ in Situation \ref{situation_finite_group}. Recall that we have a line bundle $ \cH = (H,\phi) $ on $ \cX $ which descends to an ample line bundle on the surface $ X/G. $
We first record a result for later use:
$$ \chi(\cX, \cO_\cX) = 2. $$ 
This follows from a direct computation:
\begin{enumerate}[font=\normalfont,leftmargin=*]
	\item $ h^0(\cX, \cO_\cX) = h^0(X, \cO_X) = 1 $ since $ H^0(X,\cO_X) = H^0(X,\cO_X)^G $.
	\item $ h^1(\cX, \cO_\cX) = 0 $ since $ H^1(X,\cO_X) = 0 $ and hence $ H^1(X,\cO_X)^G = 0 $.
	\item $ h^2(\cX, \cO_\cX) = h^0(\cX, \cO_\cX) = 1 $ by Serre duality for projective stacks.
\end{enumerate}
An immediate consequence is that for every line bundle $ \cL $ on $ \cX $,
\begin{equation*}
\chi(\cL,\cL) = \chi(\cX, \cL^\vee \otimes \cL) = \chi(\cX, \cO_\cX) = 2,
\end{equation*}
and hence $ \inprod{\orbv(\cL)^2}_{I\cX} = 2. $

Let's set up a few notations for later use.
\begin{notation}\label{notation_K3/G}
	Let $ \cX = [K3/G] $ in Situation \ref{situation_finite_group}. Label the conjugacy classes in $ G $ by 
	$$ [g_0], [g_1], \dots, [g_l] $$
	where $ [g_0] = [1] $ is trivial. We choose a representative element $ g_i $ in each class once and for all. For each $ g \in G, $ denote its order by $ n_g. $ 
	For every point $ x $ in the set $ X^g $ fixed by an element $ g \in G $, the linear automorphism 
	$$ dg_x: T_xX \to T_xX $$ 
	has two eigenvalues $ \lambda_{g,x} $ and $ \lambda_{g,x}^{-1} $ which are primitive $ n_g $-th roots of unity.
	Each element $ 1 \neq g \in G $ determines a finite non-empty fixed point set $ X^g $ of size
	$$ s_g = |X^g|. $$ 
	For each $ 1 \leq i \leq l, $ we have
	\begin{equation*}
	X^{g_i} = \{x_{i,1},x_{i,2},\dots,x_{i,s_{g_i}}\}.
	\end{equation*}
	Note that $ X^{g_1}, \cdots,  X^{g_l} $ are not distinct when $ l > 1 $. Indeed, for any triple $ i,j,k, $
	\begin{equation*}
	g_j = g_i^k \quad \Rightarrow \quad X^{g_i} \subset X^{g_j}. 
	\end{equation*}
	The union of these $ X^{g_i} $ consist of $ r $ non-free $ G $-orbits on $ X $, and they correspond to the orbifold points on $ \cX, $ and the singular points on $ Y $: 
	\begin{equation*}
	Gx_1, \dots, Gx_r \ \text{on} \ X \quad \leftrightarrow \quad \{p_1, \dots, p_r\} \ \text{on} \ \cX \quad \leftrightarrow \quad y_1, \dots, y_r \ \text{on} \ Y.
	\end{equation*}
	Here a point $ x_k $ is chosen in each $ G $-orbit on $ X $ once and for all. Hence, each $ x_{ij} $ is one of these $ x_k's. $
	Note that the stabilizer
	\begin{equation*}
	G_{k} = \Stab_G(x_k)
	\end{equation*}
	of $ x_k $ in $ G $ is nontrivial for all $ 1 \leq k \leq r. $
	Take an element $ g_i \neq 1. $ The finite set $ X^{g_i} $ is invariant under the centralizer $ Z_{g_i} $ of $ g_i $. The number of orbits on $ X^{g_i} $ under the $ Z_{g_i} $-action is
	\begin{equation*}
	m_i = |X^{g_i}/Z_{g_i}|,
	\end{equation*} 
	and label them by
	\begin{equation*}
	O_{i,1}, O_{i,2}, \dots, O_{i,m_i}.
	\end{equation*}
	For each $ 1 \leq j \leq m_i, $ there is a unique index $ 1 \leq k_{ij} \leq r $ such that $ x_{k_{ij}} $ is in $ X^{g_i} $ and
	\begin{equation*}
	O_{ij} = Z_{g_i} \cdot x_{k_{ij}}.
	\end{equation*}
	Denote the total number of these $ Z_{g_i} $ orbits $ O_{ij} $ on $ X^{g_i} $ by
	$$ m = m_1 + \cdots + m_l. $$ 
	For each point $ x_{k_{ij}} \in X^{g_i} $, define
	\begin{equation*}
	G_{ij} = \Stab_{Z_{g_i}}(x_{k_{ij}}) = Z_{g_i} \cap \Stab_G(x_{k_{ij}}).
	\end{equation*}
	Note that each $ G_{ij} $ is a subgroup of $ Z_{g_i} $ with size 
	$$ |G_{ij}| = |Z_{g_i}|/|O_{ij}|. $$
	The inertia stack of $ \cX $ is given by $ I\cX = \cX \coprod I_t \cX, $ where
	\begin{equation*}
	I_t \cX = \coprod_{i=1}^l \{g_i\} \times \left[ X^{g_i}/Z_{g_i} \right] = \coprod_{i=1}^l \coprod_{j=1}^{m_i} \{g_i,x_{k_{ij}}\} \times BG_{ij}.
	\end{equation*}
\end{notation}
\begin{ex}
	If $ G = \mu_2 = \{1, \sigma\}$ is generated by a Nikulin involution i.e., a symplectic involution $ \sigma, $ then there are eight fixed points of $ \sigma $,
	\begin{equation*}
	x_1, \cdots, x_8,
	\end{equation*}
	which correspond to eight $ A_1 $ singularities on $ Y. $ So $ l = 1, m_1 = 8, $ and each orbit $ O_{1,j} = \{x_j\} $ with stabilizer $ G_{1,j} = \mu_2. $ The twisted sectors are
	\begin{equation*}
	I_t \cX = \coprod_{j=1}^{8} \{\sigma,x_j\} \times B\mu_2.
	\end{equation*}
\end{ex}


\subsection{Orbifold HRR formula}\label{sec_HRR_K3/G}
In this section we work out an explicit orbifold HRR formula for $ \cX = [K3/G] $ in Notation \ref{notation_K3/G}. 
Recall that the Chern character map on the $ K3 $ surface $ X $ gives a ring isomorphism 
\begin{equation*}
\ch: N(X) \to R(X) = R^0(X) \oplus R^1(X) \oplus R^2(X) \cong \ZZ \oplus \Pic(X) \oplus \ZZ,
\end{equation*}
where $ \Pic(X) $ is a free $ \ZZ $-module with rank $ \rho(X) $ which depends on the surface $ X. $ 
The complex numerical Chow ring is then identified as
$$ R(I\cX)_\CC = R(\cX) \oplus R(I_t\cX)_\CC \cong R(X)^G \oplus \CC^m \cong \ZZ \oplus \Pic(X)^G \oplus \ZZ \oplus \CC^m $$
(Here $ R(\cX) $ is not tensored by $ \CC $.) where $ R(X)^G $ denotes the $ G $-invariant part of $ R(X) $ under the induced action of $ G $ on $ R(X) $, and $ \Pic(X)^G $ is generated by isomorphism classes of $ G $-invariant line bundles.

Now let's compute the orbifold Chern character map
\begin{equation*}
\orbch: N(\cX) \to R(X)^G \oplus \CC^m.
\end{equation*}
Take a vector bundle $ \cV = (V, \phi) $ on $ \cX $. For each point $ x_{k_{ij}}, $ the equivariant structure $ \phi $ on $ V $ restricts to a representation
\begin{equation*}
\phi_{ij}: G_{ij} \to \GL(V_{ij})
\end{equation*}
of the group $ G_{ij} $ on the fiber $ V_{ij} $ of $ V $. By the explicit formula for the orbifold Chern charcter in \cite[Section 3.3]{chen2023orbifold}, we have
\begin{equation*}
\orbch(\cV) = \left(\ch(V) , \left(\chi_{\phi_{ij}}(g_i)\right) \right)
\end{equation*}
where $ 1 \leq i \leq l, 1 \leq j \leq m_i. $
For a sheaf $ \cE = (E, \phi) $ on $ \cX, $ we resolve it by any finite locally free resolution $ \cE^{\boldsymbol{\cdot}} \to \cE \to 0 $ with each $ \cE^k = (E^k, \phi^k) $, and hence we have
\begin{equation*}
\orbch(\cE) = \sum_{k} (-1)^k \orbch(\cE^k).
\end{equation*}

Since $ G $ is a finite group, the Todd class $ \td_\cX = \td_X. $ The orbifold Todd class of the tangent complex of $ \cX $ is given by
\begin{equation*}
\orbtd_\cX = \left(\td_X,\ \bigoplus_{i,j} \frac{1}{e^{\rho_{g_i}}(T_{ij} X)}\right) = \left((1,0,2),\ \bigoplus_{i,j} \frac{1}{2-2\Real(\lambda_{ij})}\right),
\end{equation*}
where $ T_{ij} X $ is the tangent space of $ X $ at $ x_{k_{ij}} $, and $ \lambda_{ij} $ is either one of the two conjugate eigenvalues of $ g_i $ at the point $ x_{k_{ij}} $ which is also a primitive $ n_{g_i} $-th root of unity. Here each $ e^{\rho_{g_i}} $ is a twisted Euler class map. See \cite[Section 3.3]{chen2023orbifold} for the definition of the twisted Euler class map and its explicit formula. 

Consider a connected proper smooth quotient stack $ \cX $. The orbifold Euler pairing
\begin{equation*}
	\chi: N(\cX) \times N(\cX) \to \ZZ
\end{equation*}
is defined by 
\begin{equation*}
	\chi(\cE,\cF) = \sum_i (-1)^i \dim \Ext^i (\cE, \cF)
\end{equation*}
for coherent sheaves $ \cE $ and $ \cF $ on $ \cX $ and extended linearly.
In \cite[Section 3.4]{chen2023orbifold}, an orbifold Mukar vector map 
\begin{equation*}
	\orbv: N(\cX) \to R(I\cX)_{\CC}
\end{equation*}
and an orbifold Mukai pairing 
\begin{equation*}
	\inprod{\cdot \ {,}\ \cdot}_{I\cX}: R(I\cX)_\CC \times R(I\cX)_\CC \to \CC
\end{equation*}
were defined and used to derive an orbifold HRR formula
\begin{equation*}
	\chi(x,y) = \inprod{\orbv(x),\orbv(y)}_{I\cX}.
\end{equation*}

Now we work out the explicit orbifold HRR formula for our model $ \cX = [K3/G]. $

The orbifold Mukai vector map
\begin{equation*}
\orbv: N(\cX) \to R(X)^G \oplus \CC^m
\end{equation*}
is given by 
\begin{equation*}
\orbv(\cE) = \orbch(\cE) \sqrt{\td_{I\cX}} = \orbch(\cE) ((1,0,1),\underbrace{1, \dots, 1}_m ) = \left(v(E), \left(\chi_{\phi_{ij}}(g_i)\right) \right)
\end{equation*}
for a sheaf $ \cE = (E, \phi) $ on $ \cX $, where 
$$ v(E) = (\rk(E), \ch_1(E), \rk(E) + \ch_2(E)) $$ 
in $ R(X)^G $ is the Mukai vector of the sheaf $ E $ on $ X $.

The orbifold Mukai pairing is
\begin{align}\label{eq_orbv_pairing_K3/G}
\inprod{\orbv,\orbw}_{I\cX} & = \frac{\inprod{v,w}_X}{|G|} + \frac{1}{2} \sum_{i=1}^{l} \sum_{j=1}^{m_i} \frac{\overline{v}_{ij} w_{ij}}{|G_{ij}|(1-\Real(\lambda_{ij}))}
\end{align}
for all $ \orbv = (v,(v_{ij})) $ and $ \orbw = (w,(w_{ij})) $ in $ R(I\cX)_\CC \cong R(X)^G \oplus \CC^m $, where 
$$ \inprod{v,w}_X = \int_X v^\vee w = v_0 w_2 - v_1 w_1 + v_2 w_0 $$ in $ \ZZ $
is the integral Mukai pairing of $ v $ and $ w $ in $ R(X). $

For every pair of sheaves $ \cE = (E, \phi) $ and $ \cF = (F, \psi) $ on $ \cX, $ the orbifold HRR formula is
\begin{equation*}
\chi(\cE,\cF) = \inprod{\orbv(\cE),\orbv(\cF)}_{I\cX} = \frac{\inprod{v(E),v(F)}_X}{|G|} + \frac{1}{2} \sum_{i=1}^l \sum_{j=1}^{m_i} \frac{\chi_{\phi_{ij}}(g_i^{-1}) \chi_{\psi_{ij}}(g_i)}{|G_{ij}|(1-\Real(\lambda_{ij}))};
\end{equation*}
in particular, if $ \cE = \cF, $ then
\begin{equation}\label{eq_orbv_squared}
\chi(\cE,\cE) = \inprod{\orbv(\cE)^2}_{I\cX} = \frac{\inprod{v(E)^2}_X}{|G|} + \frac{1}{2} \sum_{i=1}^l \sum_{j=1}^{m_i} \frac{|\chi_{\phi_{ij}}(g_i)|^2}{|G_{ij}|(1-\Real(\lambda_{ij}))}.
\end{equation}
\begin{rmk}
	Consider any coherent sheaf $ \cE $ on $ \cX. $ Since $ \chi(\cE,\cE) $ is an integer, the double summation in (\ref{eq_orbv_squared}) must be a rational number although each term is irrational in general. Moreover, $ \chi(\cE,\cE) $ is always an even integer as we will see shortly.
\end{rmk}
We first generalize Proposition \ref{prop_K3_smoothness_and_dim} for $ K3 $ surfaces to $ [K3/G]. $
\begin{prop}\label{prop_K3/G_smoothness_dim}
	Let $ \cX = [K3/G] $ in Situation \ref{situation_finite_group}. Let $ x $ be an element in $ N(\cX) $ with orbifold Mukai vector $ \orbv $ in $ R(I\cX)_\CC $. Then $ M^s(\cX,x) $ is either empty or a smooth quasi-projective scheme with
	\begin{equation*}
	\dim M^s(\cX,x) = 2 - \inprod{\orbv^2}_{I\cX}.
	\end{equation*}
	Denote the Mukai vector of $ y = p^N x $ in $ N(X) $ by $ v = (r, c_1, s) $ in $ R(X)^G $ with $ d = \deg(y). $ Suppose the pair $ (y,H) $ satisfies either of the following conditions:
	\begin{enumerate}[font=\normalfont,leftmargin=2em]
		\item $ \gcd(r, d, s) = 1 $.
		\item $ y $ is primitive and $ H $ is $ y $-generic.
	\end{enumerate}
	Then $ M(\cX,x) = M^s(\cX,x). $ If $ M(\cX,x) $ is not empty, then it is a smooth projective scheme of dimension $ 2 - \inprod{\orbv^2}_{I\cX}. $
\end{prop}
\begin{proof}
	Suppose $ M_\cH^s(\cX,x) $ is not empty. Then it is a quasi-projective scheme by construction. Take a stable sheaf $ \cE $ on $ \cX $ with numerical class $ x. $ By a similar argument to that in \cite[Proposition 1.2.7]{huybrechts2010geometry}, an endomorphism $ \cE \to \cE $ is either zero or invertible. Since the automorphism group of a stable sheaf is $ \CC^*, $ we must have
	\begin{equation*}
	\End(\cE) \cong \CC.
	\end{equation*}
	In particular, there is an isomorphism 
	$$ H^0(\cX, \cO_\cX) \cong \CC \xrightarrow{\sim} \Ext^0(\cE, \cE) $$
	of $ \CC $-vector spaces.
	Since the canonical sheaf $ \omega_\cX \cong \cO_\cX, $ Serre duality for smooth projective stacks gives an isomorphism
	\begin{equation*}
	\Ext^2(\cE, \cE) \xrightarrow{\sim} H^2(\cX, \cO_\cX).
	\end{equation*}
	The Picard scheme $ \Pic(\cX) = \Pic^G(X) $ is a group scheme, which is always smooth in characteristic zero. Therefore, by Proposition \ref{prop_smoothness_criterion}, $ M^s(\cX,x) $ is smooth at the closed point $ t $ corresponding to $ \cE. $ So we can compute its dimension at the tangent space $ T_tM^s(\cX,x) \cong \Ext^1(\cE,\cE) $. Recall the orbifold Euler pairing
	\begin{align*}
	\chi(\cE, \cE) = \sum_{i=0}^2 \dim \Ext^i(\cE, \cE) = 2 - \dim \Ext^1(\cE, \cE).
	\end{align*}
	By the orbifold HRR formula, we have
	\begin{equation*}
	\dim \Ext^1(\cE, \cE) = 2 - \chi(\cE, \cE) = 2 - \inprod{\orbv(\cE)^2}_{I\cX}.
	\end{equation*}
	Now suppose either of the two conditions in Proposition \ref{prop_K3/G_smoothness_dim} is satisfied. If $ M(\cX,x) $ is empty, there is nothing to prove. Suppose it's not empty. Take any semistable sheaf $ \cF = (F, \psi) \in M(\cX,x) $. By Proposition \ref{prop_Gieseker_relation}, the coherent sheaf $ F $ on $ X $ is also semistable. By Proposition \ref{prop_K3_smoothness_and_dim}, $ F $ is also stable, and hence $ \cE $ is stable as well again by Proposition \ref{prop_Gieseker_relation}. The last statement follows since $ M(\cX,x) $ is a projective scheme by construction.
\end{proof}
\begin{rmk}
	The results in Proposition \ref{prop_K3/G_smoothness_dim} do not say that the moduli space $ M^s(\cX,x) $ is irreducible, which will be proved later in Theorem \ref{thm_K3/G} under a slightly stronger assumption, i.e., $ r >0, $ and either $ d > 0 $ or $ \gcd(r,d) = 1. $
\end{rmk}
Let's apply Proposition \ref{prop_K3/G_smoothness_dim} to a few examples. 
\begin{ex}
	Let $ \cL $ be a line bundle on $ \cX $ with numerical class $ \gamma(\cL) $ in $ N(\cX). $ Since $ \inprod{\orbv(\cL)^2} = 2, $ Proposition \ref{prop_K3/G_smoothness_dim} implies
	\begin{equation*}
	\dim  M(\cX, \gamma(\cL)) = 0,
	\end{equation*}
	since there are no strictly semistable sheaves of rank one. Indeed, each moduli space $ M(\cX, \gamma(\cL)) $ is a point in the smooth Picard scheme $ \Pic(\cX) $ of $ \cX $.
\end{ex}
\begin{ex}
	Let $ \cO_p $ be the structure sheaf of a generic point $ p $ on $ \cX $ corresponding to a free orbit $ Gx $ on $ X. $ Then $ \cO_p $ is a stable sheaf on $ \cX $ with an orbifold Mukai vector
	\begin{equation*}
	\orbv(\cO_p) = (0,0,|G|,\underbrace{0, \dots, 0}_m)
	\end{equation*}
	in $ R(I\cX)_\CC. $ Therefore, $ \inprod{\orbv(\cO_p)^2} = 0, $ which implies
	\begin{equation*}
	\dim  M(\cX, \gamma(\cO_p)) = 2.
	\end{equation*}
	Indeed, we can identify
	$$ M(\cX, \gamma(\cO_p)) \cong \GHilb X, $$ the $ G $-Hilbert scheme of free orbits on $ X. $ Hence $ M = M(\cX, \gamma(\cO_p)) $ is a $ K3 $ surface which is the minimal resolution of the surface $ X/G. $ Moreover, $ M $ is a fine moduli space of stable sheaves on $ \cX $.
\end{ex}

The following result is immediate from the orbifold HRR formula for $ [K3/G] $.
\begin{prop}
	Let $ \cX = [K3/G] $ in Situation \ref{situation_finite_group} with Notation \ref{notation_K3/G}. Then
	\begin{equation}\label{eq_K3/G_identity}
	\frac{1}{|G|} + \frac{1}{4} \sum_{i=1}^l \sum_{j=1}^{m_i} \frac{1}{|G_{ij}|(1-\Real(\lambda_{ij}))} = 1.
	\end{equation}
\end{prop}
\begin{proof}
	Consider $ \cE = \cO_\cX $ in (\ref{eq_orbv_squared}). The result follows from the fact $ \chi(\cX,\cO_\cX) = 2. $ 
\end{proof}

A result of Nikulin \cite{nikulin1979finite} says that if a $ K3 $ surface admits a symplectic automorphism of finite order $ n $, then $ n \leq 8. $ In \cite[Section 1]{mukai1988finite}, Mukai proved that the number of points $ f_n $ fixed by a symplectic automorphism only depends on its order $ n $, and computed all such numbers in Table \ref{table_number_fixed_points} on page \pageref{table_number_fixed_points}.

We now use equation (\ref{eq_K3/G_identity}) to reproduce these numbers $ f_n $. Note that $ f_n = |X^G|, $ where $ G $ is the cyclic group generated by a symplectic automorphism of order $ n. $
We need not assume that $ |X^G| $ only depends on the order of $ G $, which is rather a result from our proof of the following

\begin{cor}\label{cor_table_fixed_points}
	Let $ \cX = [K3/G] $ in Situation \ref{situation_finite_group} where $ G \cong \ZZ/n\ZZ $. If $ n $ is a prime $ p $, then the number of fixed points of $ G $ is 
	$$ |X^G| = \frac{24}{p+1}. $$ 
	If $ n = 4, 6, $ and $ 8, $ then $ |X^G| = 4, 2 $ and $ 2 $ respectively. Thus we reproduce Table \ref{table_number_fixed_points} on page \pageref{table_number_fixed_points}. 
\end{cor}
\begin{proof}
	Let $ g $ be a generator of $ G, $ i.e., 
	$$ G = \{1, g, g^2, \dots, g^{n-1}\}. $$
	Then $ X^G = X^g $. 
	For each $ 1 \leq i \leq n-1, $ let $ s_i $ denote the number of fixed points of $ g^i. $ 
	Suppose $ n $ is a prime. Then every nontrivial $ g^i $ has the same fixed point set $ X^G $, and hence $ s_i = s_1 = |X^G| $ for all $ 1 \leq i \leq n-1. $ Since $ G $ is abelian, each centralizer $ Z_{g^i} = G $, so each stack 
	$$ [X^{g^i}/Z_{g^i}] = [X^G/G]. $$ 
	Therefore, in equation (\ref{eq_K3/G_identity}), the number of orbits $ m_i = |X^G| $ for all $ 1 \leq i \leq n-1 $, and each stabilizer $ G_{ij} $ has size $ |G_{ij}| = n. $
	Note that for each $ 1 \leq k \leq n-1, $
	\begin{equation*}
	\lambda_{k} = \exp(2k\pi i/n)
	\end{equation*}
	is an eigenvalue of $ g^k. $ 
	Equation (\ref{eq_K3/G_identity}) then becomes
	\begin{equation}\label{eq_K3/G_identity_prime_order}
	\frac{1}{n} + \frac{|X^G|}{4n} \sum_{k=1}^{n-1} \frac{1}{1-\cos(2k\pi /n)} = 1.
	\end{equation}
	Recall a trigonometric identity
	\begin{equation}\label{eq_trig_identity}
	\sum_{k=1}^{n-1} \frac{1}{1-\cos(2k\pi /n)} = \frac{n^2-1}{6}.
	\end{equation}
	Plugging (\ref{eq_trig_identity}) into (\ref{eq_K3/G_identity_prime_order}) yields
	\begin{equation*}
	|X^G| = \frac{24}{n+1}.
	\end{equation*}
	This gives $ f_2,f_3,f_5,f_7 $ in Table \ref{table_number_fixed_points}. The remaining three cases can be computed directly. For a non-prime $ n, $ we use the same argument, but now for each $ 1 \leq i \leq n-1, $ we have
	\begin{equation*}
	\sum_{j=1}^{m_i} \frac{1}{|G_{ij}|} = \sum_{j=1}^{m_i} \frac{|O_{ij}|}{|G|} = \frac{s_i}{n}.
	\end{equation*}
	Therefore, equation (\ref{eq_K3/G_identity}) becomes
	\begin{equation}\label{eq_K3/G_identity_non_prime_order}
	\frac{1}{n} + \frac{1}{4n} \sum_{k=1}^{n-1} \frac{s_k}{1-\cos(2k\pi /n)} = 1.
	\end{equation}
	Note that for each $ 1 \leq k \leq n-1, $ the element $ g^k $ has order 
	$$ n_k = \frac{n}{\gcd(n,k)}, $$ 
	so the number of fixed points $ s_k = f_{n_k}$ when $ n_k $ is a prime. Now applying (\ref{eq_K3/G_identity_non_prime_order}) to $ n = 4,$ we get
	\begin{align*}
	\frac{1}{4} + \frac{1}{16} \left(|X^G| + \frac{f_2}{2} + |X^G|\right) = 1,
	\end{align*}
	which gives $ |X^G| = 4. $ This implies all groups $ G \cong \ZZ/4\ZZ $ acting faithfully and symplectically on $ X $ have the same number of fixed points, and gives $ f_4 = 4 $ in Table \ref{table_number_fixed_points}.
	Carrying out the same procedure for $ n = 6 $ and $ 8, $ we have
	\begin{align*}
	& \frac{1}{6} + \frac{1}{24} \left(\frac{f_6}{1/2} + \frac{f_3}{3/2} + \frac{f_2}{2} + \frac{f_3}{3/2} + \frac{f_6}{1/2}\right) = 1, \ \text{and}\\
	& \frac{1}{8} + \frac{1}{32} \left(\frac{f_8}{1-1/\sqrt{2}} + f_4 + \frac{f_8}{1+1/\sqrt{2}} + \frac{f_2}{2} + \frac{f_8}{1+1/\sqrt{2}} + f_4 + \frac{f_8}{1-1/\sqrt{2}}\right) = 1,
	\end{align*}
	which gives $ f_6 = f_8 = 2, $ thus completing Table \ref{table_number_fixed_points} for $ 2 \leq n \leq 8 $. 
\end{proof}
\vskip 1pt

\subsection{Derived McKay correspondence}
Let $ p $ denote a generic point on $ \cX = [K3/G] $. Consider again the moduli space 
$$ M = M(\cX,\gamma(\cO_p)). $$
We have seen that $ M = \GHilb X, $ which is a fine moduli space. Therefore, the universal family on $ \cX \times M $ induces a Fourier-Mukai transform
\begin{equation*}
\Phi: \D(\cX) \xrightarrow{\sim} \D(M)
\end{equation*}
which is an equivalence between derived categories known as the derived McKay correspondence in \cite{bridgeland2001mckay}. Passing to the numerical Grothendieck rings and the numerical Chow rings, $ \Phi $ induces a commutative diagram
\vspace{1pt}
\begin{equation*}\label{diag_Fourier_Mukai}
\begin{tikzcd}[column sep=2.5em,row sep=3em,every label/.append style={font=\normalsize}]
N(\cX) \arrow[r,"\Phi^N"{outer sep = 2pt}] \arrow[d, "\orbv"{left, outer sep = 2pt}]
& N(M) \arrow[d, "v"{right, outer sep = 2pt}] \\
R(I\cX)_\CC \arrow[r,"\Phi^R"{outer sep = 2pt}] & R(M)_\CC,
\end{tikzcd}\vspace{3pt}
\end{equation*}
where $ \Phi^N $ and $ \Phi^R $ are isomorphisms of rings and $ \CC $-algebra respectively, $ \orbv $ and $ v $ are injections of abelian groups, and each arrow is compatible with the canonical pairings on the sources and targets. For example, $ \Phi^N $ preserves orbifold Euler pairings:
\begin{equation*}
\chi(x, y) = \chi(\Phi^N(x),\Phi^N(y))
\end{equation*}
in $ \ZZ $ for $ x $ and $ y $ in $ N(\cX), $ and $ \Phi^R $ preserves orbifold Mukai pairings:
\begin{equation*}
\inprod{\orbv,\orbw}_{I\cX} = \inprod{\Phi^R(\orbv),\Phi^R(\orbw)}_M
\end{equation*}
in $ \CC $ for $ \orbv $ and $ \orbw $ in $ R(I\cX)_\CC. $ This implies the following
\begin{prop}
	Let $ \cX = [K3/G] $ in Situation \ref{situation_finite_group}. Let $ \cE $ be a coherent sheaf on $ \cX $. Then the integer $ \inprod{\orbv(\cE)^2}_{I\cX} $ is even.
\end{prop}
\begin{proof}
	Consider the minimal resolution $ M \to X/G. $ We have 
	$$ \inprod{\orbv(\cE)^2}_{I\cX} = \inprod{\Phi^R(\orbv(\cE))^2}_M = \inprod{v(\Phi^N(\cE))^2}_M, $$ 
	which is even since the self-intersection on $ R^1(Y) $ is even for any $ K3 $ surface $ Y $.  
\end{proof}
\vspace{1pt}

\subsection{Bridgeland stability conditions}\label{sec_Bridgeland}
In this section we review Bridgeland stability conditions on triangulated categories, their constructions on a $ K3 $ surface $ X $, and the induced stability conditions on the derived category of the quotient stack $ \cX = [X/G] $ where $ G $ is a finite group. We will relate Gieseker moduli spaces of sheaves on $ \cX $ to Bridgeland moduli spaces of complexes in $ \D(\cX). $ This connection is a key ingredient in our proof of the main theorem later.

We first review the Bridgeland stability conditions for triangulated categories introduced in \cite{bridgeland2007stability}. Fix a triangulated category $ \cD $. Let $ K(\cD) $ denote the Grothendieck group of $ \cD. $ 

\begin{defn}[{\cite[Definition 5.1]{bridgeland2007stability}}]
	A stability condition $ \sigma = (Z, \cP) $ on $ \cD $ consists of an additive map $ Z: K(\cD) \to \CC $ called the \tb{central charge} of $ \sigma $, and a slicing $ \cP $ of $ \cD $ such that for each \tb{phase} $ \theta $ in $ \RR, $ if a nonzero object $ E $ is in $ \cP(\theta), $ then
	\begin{equation*}
	Z(E) = m(E) e^{i\pi\theta}
	\end{equation*}
	for some $ m(E) > 0 $ called the \tb{mass} of $ E. $
\end{defn}
If there exists a stability condition on $ \cD $, then there are notions of stability and semistability for the objects in $ \cD. $
\begin{defn}
	Let $ \sigma = (Z, \cP) $ be a \tb{stability condition} on $ \cD. $ For each $ \theta $ in $ \RR, $ the nonzero objects of $ \cP(\theta) $ are said to be \tb{$ \sigma $-semistable} of phase $ \theta $, and the simple\footnote{A simple object in a category is an object that has no nonzero proper subobjects. They are also called minimal objects in the literature.} objects of $ \cP(\theta) $ are said to be \tb{$ \sigma $-stable}.
\end{defn}
Alternatively, one can define stability conditions via a stability function on the \tb{heart} $ \cA $ of a \tb{bounded $ t $-structure} on $ \cD $, which is an abelian subcategory of $ \cD. $
\begin{defn}[{\cite[Definition 2.1]{bridgeland2007stability}}]
	A \tb{stability function} on an abelian category $ \cA $ is an additive map $ Z: K(\cA) \to \CC $ such that for all nonzero objects $ E $ in $ \cA, $ the complex number $ Z(E) $ lies in the union of the strict upper half plane and the negative real axis, i.e.,
	\begin{equation*}
	Z(E) = m(E) e^{i\pi\theta(E)}
	\end{equation*}
	where $ m(E) > 0 $ and $ \theta(E) \in (0,1]. $
\end{defn}
\begin{defn}
	Let $ Z: K(\cA) \to \CC $ be a stability function on an abelian category $ \cA $. A nonzero objects $ E $ of $ \cA $ is said to be (semi)stable with respect to $ Z $ if
	\begin{equation*}
	\theta(F) \ (\leq)\ \theta(E)
	\end{equation*} 
	for all nonzero proper subobjects $ F \subset E. $ A \tb{Harder-Narasimhan (HN) filtration} of a nonzero object $ E $ in $ \cA $ is a chain of subobjects
	\begin{equation*}
	0 = E_0 \subset E_1 \subset \cdots \subset E_n = E
	\end{equation*}
	such that each quotient $ F_i = E_i/E_{i-1} $ is a semistable object in $ \cA $ with
	\begin{equation*}
	\theta(F_1) > \theta(F_2) > \cdots > \theta(F_n).
	\end{equation*}
	The stability function $ Z $ is said to have the \tb{HN property} if every nonzero object of $ \cA $ has an HN filtration.
\end{defn}
If $ \sigma = (Z, \cP) $ is a stability condition on $ \cD $, then the abelian subcategory $ \cA = \cP((0,1]) $ of $ \cD $ is the heart of the $ t $-structure $ \cP(>0), $ and the central charge $ Z: K(\cD) \to \CC $ gives a stability function $ Z: K(\cA) \to \CC. $ This yields two equivalent ways to impose stability conditions on $ \cD $.
\begin{prop}[{\cite[Proposition 5.3]{bridgeland2007stability}}]
	A stability condition $ \sigma = (Z, \cP) $ on $ \cD $ is the same as a pair $ (\cA, Z) $ where $ \cA $ is the heart of a bounded $ t $-structure on $ \cD $ and $ Z $ is a stability function on $ \cA $ with the HN property.
\end{prop}

For be a smooth projective variety $ X $ over $ \CC, $ it's a non-trivial question whether there exist stability conditions on the derived categories $ \D(X) $ of $ X $. The answer is positive when $ X $ is a $ K3 $ surface.

\begin{constrn}[Bridgeland stability conditions on $ K3 $ surfaces]
	Consider a $ K3 $ surface $ X $. In \cite{bridgeland2008stability}, Bridgeland constructed stability conditions on $ X $ as follows. Fix a pair $ (\beta, \omega) $ in $ R^1(X)_\RR $ such that $ \omega \in \Amp(X). $ Let $ (\beta \omega) = \deg(\beta \omega)$ denote the degree of the intersection of $ \beta $ and $ \omega $ on $ R^1(X)_\RR $.
	Every torsion-free sheaf $ E $ on $ X $ has an HN filtration
	\begin{equation*}
	0 = E_0 \subset E_1 \subset \cdots \subset E_n = E
	\end{equation*}
	such that each HN factor $ F_i = E_i/E_{i-1} $ is a $ \mu_\omega $-semistable torsion-free sheaf on $ X $ with
	\begin{equation*}
	\mu_\omega(F_1) > \mu_\omega(F_2) > \cdots > \mu_\omega(F_n).
	\end{equation*}
	Define two full additive subcategories of the abelian category $ \Coh(X) $ (also known as a \tb{torsion pair}) as follows:
	\begin{align*}
	\cT(\beta,\omega) & = \{T \in \Coh(X): \text{the torsion-free part of}\ T\ \text{has HN factors with}\ \mu_\omega > (\beta \omega) \}. \\
	\cF(\beta,\omega) & = \{F \in \Coh(X): F\ \text{is torsion-free and has HN factors with}\ \mu_\omega \leq (\beta \omega) \}.
	\end{align*}
	Define an abelian subcategory of $ \D(X) $ by
	\begin{equation*}
	\cA(\beta,\omega) = \{E \in \D(X): H^i(E) = 0 \ \text{for}\ i \neq -1, 0, H^{-1}(E) \in \cF(\beta,\omega), H^0(E) \in \cT(\beta,\omega) \}.
	\end{equation*}
	Then $ \cA(\beta,\omega) $ is the heart of a bounded $ t $-structure on $ \D(X). $ Note that $ \cT(\beta,\omega) $ is a subcategory of $ \cA(\beta,\omega) $ and contains all torsion sheaves on $ X. $
	Define an additive map
	\begin{equation*}
	Z_{\beta,\omega}: N(X) \to \CC, \quad x \mapsto -\inprod{\exp(\beta-i\omega),v(x)}
	\end{equation*}
	where $ \inprod{\cdot\ {,}\ \cdot} $ denotes the sesquilinear Mukai pairing on the complex numerical Chow ring $ R(X)_\CC. $
\end{constrn}
\begin{prop}[{\cite[Lemma 6.2]{bridgeland2008stability}}]\label{prop_Bridgeland_K3}
	Let $ X $ be a $ K3 $ surface. For a pair $ (\beta, \omega) $ in $ R^1(X)_\RR $ such that $ \omega \in \Amp(X) $ with $ \deg(\omega^2) > 2, $ the map $ Z_{\beta,\omega}: N(X) \to \CC $ is a stability function on the abelian category $ \cA(\beta,\omega) $ and hence the pair $ (Z_{\beta,\omega}, \cA(\beta,\omega)) $ gives a stability condition $ \sigma_{\beta,\omega} $ on $ \D(X). $
\end{prop}
\begin{rmk}
	Let $ X $ be a smooth projective variety over $ \CC. $ The space of all stability conditions on $ X $ is denoted by $ \Stab(X). $ The map 
	$$ \Stab(X) \to \Hom(N(X), \CC), \quad (\cA, Z) \mapsto Z $$
	is a local homeomorphism which endows each connected component of $ \Stab(X) $ a structure of a complex manifold of dimension $ \rk(N(X)) $. When $ X $ is a $ K3 $ surface, there is a \tb{distinguished component} which contains all stability conditions $ \sigma_{\beta,\omega} $ in Proposition \ref{prop_Bridgeland_K3}, and is denoted by $ \Stab^\dagger(X). $
\end{rmk}
Given a stability condition $ \sigma $ in $ \Stab(X) $ and an element $ x $ in $ N(X), $ the moduli stack $ \cM_\sigma(\D(X),x) $ has a good moduli space $ M_\sigma(\D(X),x) $ parametrizing $ S $-equivalent classes of $ \sigma $-semistable complexes in $ \D(X) $ with numerical class $ x $, which becomes a coarse moduli space if $ \cM_\sigma(\D(X),x) $ coincides with the moduli stack $ \cM_\sigma^s(\D(X),x) $ of $ \sigma $-stable complexes in $ \D(X) $ with numerical class $ x $. 
\begin{thm}[{\cite[Theorem 21.24]{bayer2021stability}}]
	Let $ X $ be a smooth projective variety $ X $ over $ \CC $. Let $ \sigma $ be a stability condition on $ \D(X) $, and let $ x $ be an element in $ N(X). $ Then $ \cM_\sigma(\D(X),x) $ is an Artin stack of finite type over $ \CC $ with a good moduli space $ M_\sigma(\D(X),x) $ which is a proper algebraic space. If $ \cM_\sigma(\D(X),x) = \cM_\sigma^s(\D(X),x), $ then it is a $ \CC^* $-gerbe over its coarse moduli space $ M_\sigma(\D(X),x) = M_\sigma^s(\D(X),x). $
\end{thm}
Now we consider a $ K3 $ surface $ X. $ Under certain conditions, the moduli space $ M_\sigma(\D(X),x) $ is a projective variety deformation equivalent to a Hilbert scheme of points on a $ K3 $ surface. Fix an element $ x $ in $ N(X). $ Then it determines a set of real codimension-one submanifolds with boundaries in the manifold $ \Stab(X) $ known as \tb{walls} as shown in \cite[Proposition 2.3]{bayer2014projectivity}.
\begin{defn-prop}[{\cite[Definition 2.4]{bayer2014projectivity}}]
	Let $ X $ be a $ K3 $ surface. Let $ x $ be an element in $ N(X). $ A stability condition $ \sigma $ in $ \Stab(X) $ is said to be \tb{$ x $-generic} if it does not lie on any of the walls determined by $ x. $ 
\end{defn-prop}
The following result is a generalization of Theorem \ref{thm_yoshioka} to Bridgeland moduli spaces of stable objects in $ \D(X) $.
\begin{thm}[{\cite[Theorem 1.1]{bottini2021stable}}]\label{thm_bottini}
	Let $ X $ be a $ K3 $ surface. Let $ x $ be an element in $ N(X) $ with a primitive Mukai vector $ v $ in $ R(X). $ Let $ \sigma $ be a stability condition in $ \Stab^\dagger(X) $ such that it is $ x $-generic. Then $ M_\sigma(\D(X),x) =  M_\sigma^s(\D(X),x) $ and 
	it is non-empty if and only if $ \inprod{v^2} \leq 2. $ In this case, $ M_\sigma(\D(X),x) $ is an irreducible symplectic manifold of dimension $ 2 - \inprod{v^2} $ deformation equivalent to a Hilbert scheme of points on a $ K3 $ surface.
\end{thm}
Under some conditions, a Bridgeland moduli space is a Gieseker moduli space.
\begin{prop}\label{prop_Bridgeland_vs_Gieseker}
	Let $ (X, H) $ be a polarized $ K3 $ surface. Let $ x $ be an element in $ N(X) $ with Mukai vector $ v = (r, c_1, s) $ in $ R(X) $ and $ d = \deg(x) $. Suppose $ r > 0, $ and either of the following conditions is satisfied:
	\begin{enumerate}[font=\normalfont,leftmargin=2em]
		\item $ d > 0. $
		\item $ \gcd(r, d) = 1. $
		\item $ (r, c_1) $ is primitive and $ H $ is $ x $-generic.
	\end{enumerate}
	Then there exists an integer $ b $ such that for all real numbers $ t \gg 0 $, there is a stability condition $ \sigma_t = \sigma_{bh,th} $ in the distinguished component $ \Stab^\dagger(X) $ such that there is an isomorphism
	\begin{equation*}
	M_H^{(s)}(X,x) \cong M_{\sigma_t}^{(s)}(\D(X), x)
	\end{equation*}
	between Gieseker and Bridgeland moduli spaces on $ X $.
	In case $ (2) $ and $ (3) $, we also have $ M_H(X,x) = M_H^s(X,x). $
\end{prop}
\begin{proof}
	Let $ h = c_1(H) $ in $ R^1(X) $. Since the self-intersection on $ R^1(X) $ is even and $ H $ is ample, we have 
	$$ \deg(H) = \deg(h^2) \geq 2. $$ Since $ r > 0, $ the element $ x $ has a slope 
	$$ \mu(x) = d/r \in \QQ. $$ Choose an integer $ b < \mu(x)/\deg(H). $ Let $ B = H^{\otimes b} $ with first Chern class $ \beta = c_1(B) $ in $ R^1(X). $ Then we have
	$$ \deg(B) = \deg(H) b < \mu(x). $$ 
	Since $ \deg(H) \geq 2, $ by Proposition \ref{prop_Bridgeland_K3}, there is a stability condition 
	$$ \sigma_t = \sigma_{bh,th} \in \Stab^\dagger(X) $$
	for all $ t > 1. $  
	Under the conditions $ \deg(B) < \mu(x) $ and $ r > 0, $ by \cite[Exercise 6.27]{macri2019lectures}, we can choose $ t_0 > 1 $ such that for all $ t \geq t_0, $ a complex $ E $ in $ \D(X) $ with numerical class $ x $ is $ \sigma_t $-(semi)stable if and only if $ E $ is a $ \beta $-twisted Gieseker $ H $-(semi)stable sheaf on $ X $ as in \cite[Definition 14.1]{bridgeland2008stability}, which is equivalent to the condition that $ E \otimes B^\vee $ on $ X $ is Gieseker $ H $-(semi)stable. Therefore, we have
	\begin{equation*}
	M_{\sigma_t}^{(s)}(\D(X), x) \cong M_H^{(s)}(X,x \cdot \gamma(B)^{-1}).
	\end{equation*}
	If $ d > 0, $ then the integer $ b $ can be chosen as zero and we are done. Suppose either $ \gcd(r, d) = 1, $ or $ (r, c_1) $ is primitive and $ H $ is $ x $-generic. By Lemma \ref{lem_slope_Gieseker} and \ref{lem_slope_semistable_implies_stable}, slope semistability is the same as Gieseker semistability. Since tensoring by the line bundle $ B $ preserves slope stability, this implies
	\begin{equation*}
	M_H^{(s)}(X,x) = M_H^{(s)}(X,x \cdot \gamma(B)^{-1}).
	\end{equation*}
	This completes the proof.
\end{proof}

\begin{constrn}[Induced Bridgeland stability conditions on quotient stacks]\label{constrn_induced_stability}
	Consider a smooth projective variety $ X $ over $ \CC $ under an action of a finite group $ G $. Then we have a quotient stack $ \cX = [X/G] $ with its derived category 
	$$ \D(\cX) = \D(\Coh(\cX)) \cong \D(\Coh^G(X)). $$ 
	The canonical morphism $ p: X \to \cX $ is smooth, and hence induces an exact functor 
	\begin{equation*}
	p^*: \Coh(\cX) \to \Coh(X)
	\end{equation*}
	and a ring homomorphism
	\begin{equation*}
	p^N: N(\cX) \to N(X), \quad \gamma(\cE) = \gamma(E,\phi) \mapsto \gamma(E)
	\end{equation*}
	for a $ G $-equivariant sheaf $ \cE = (E,\phi) $ on $ X. $
	The $ G $-action on $ X $ induces a $ G $-action on the manifold $ \Stab(X) $ via the auto-equivalence 
	$$ g^*: \D(X) \to \D(X) $$ 
	for each $ g $ in $ G. $
	Fix a $ G $-invariant stability condition $ \sigma = (Z, \cP) $
	on $ \D(X). $ By \cite[Section 2]{macri2009inducing}, it induces a stability condition $ \wtilde{\sigma} = (\wtilde{Z}, \wtilde{\cP}) $
	on $ \D(\cX) $ where the central charge $ \wtilde{Z} $ is a composition
	\begin{equation*}
	\wtilde{Z}: N(\cX) \xrightarrow{p^N} N(X) \xrightarrow{Z} \CC,
	\end{equation*}
	and the slicing $ \wtilde{\cP} $ is given by
	\begin{equation*}
	\wtilde{\cP}(\theta) = \{ \cE \in \D(\cX): p^*\cE \in \cP(\theta)\}
	\end{equation*}
	for all $ \theta \in \RR. $ Moreover, by \cite[Theorem 1.1]{macri2009inducing}, the assignment $ \sigma \mapsto \wtilde{\sigma} $ gives a closed embedding of complex manifolds 
	\begin{equation*}
	\Stab(X)^G \into \Stab(\cX).
	\end{equation*}
\end{constrn}
Let $ \cX = [X/G] $ where $ X $ is a smooth projective variety over $ \CC $ and $ G $ is a finite group acting on $ X $. Let $ p: X \to \cX $ denote the canonical morphism. Under the induced stability conditions, semistable objects in $ \D(\cX) $ pull back to semistable objects in $ \D(X). $
\begin{prop}\label{prop_Bridgeland_relation}
	Let $ \cE $ be an object in $ \D(\cX) $ with $ E = p^*\cE $ in $ \D(X) $. Then $ \cE $ is $ \wtilde{\sigma} $-semistable if and only if $ E $ is $ \sigma $-semistable. If $ E $ is $ \sigma $-stable, then $ \cE $ is $ \wtilde{\sigma} $-stable.
\end{prop}
\begin{proof}
	This is a tautology by Construction \ref{constrn_induced_stability}.
\end{proof}
Suppose $ X $ is a $ K3 $ surface. If $ \sigma $ is $ G $-invariant stability condition in $ \Stab^\dagger(X), $ then the moduli stack $ \cM_{\wtilde{\sigma}}(\D(\cX), x) $ of $ \wtilde{\sigma} $-semistable complexes in $ \D(\cX) $ with numerical class $ x $ in $ N(\cX) $ is an Artin stack and has a good moduli space $ M_{\wtilde{\sigma}}(\D(\cX), x) $.
\begin{thm}
	Let $ \cX = [X/G] $ where $ X $ is a $ K3 $ surface and $ G $ is a finite group acting on $ X $. Let $ \sigma $ be a $ G $-invariant stability condition in $ \Stab^\dagger(X). $ For every element $ x $ in $ N(\cX) $, the moduli stack $ \cM_{\wtilde{\sigma}}(\D(\cX), x) $ is an Artin stack of finite type over $ \CC $ with a proper good moduli space.
\end{thm}
\begin{proof}
	This is {\cite[Theorem 3.22]{beckmann2020equivariant}} where $ X $ is a $ K3 $ surface.
\end{proof}
\begin{rmk}
	If there are no strictly $ \wtilde{\sigma} $-semistable complexes in $ \D(\cX) $ with numerical class $ x $, then the good moduli space $ M_{\wtilde{\sigma}}(\D(\cX), x) $ is also a coarse moduli space.
\end{rmk}

\subsection{Proof of the main theorem}\label{sec_proof_main_thm}
In this section we prove the main theorem, i.e., Theorem \ref{thm_main} stated in the introduction. 

Let $ \cX = [K3/G] $ in Situation \ref{situation_finite_group}. Recall that we have a line bundle $ \cH $ on $ \cX $ which descends to an ample line bundle on the projective surface $ X/G. $ The natural morphism $ p: X \to \cX $ pulls back $ \cH $ to a $ G $-invariant ample line bundle $ H $ on $ X. $ We have a commutative diagram 
\vspace{3pt}
\begin{equation*}
\begin{tikzcd}[column sep=2.5em,row sep=3em,every label/.append style={font=\normalsize}]
N(\cX) \arrow[r, "p^N"{outer sep = 2pt}] \arrow[d, "\orbv"{left, outer sep = 2pt}] 
& N(X) \arrow[d, "v"{right, outer sep = 2pt}, ] \\ 
R(I\cX)_\CC \arrow[r, "p^R"{outer sep = 2pt}] & R(X),
\end{tikzcd}\vspace{6pt}
\end{equation*}
where the pullbacks $ p^N $ and $ p^R $ are ring homomorphisms, the Mukai vector map $ v $ and the orbifold Mukai vector map $ \orbv $ are additive maps. Recall that
$$ R(I\cX)_\CC = R(\cX) \oplus R(I_t\cX)_\CC \cong R(X)^G \oplus \CC^m. $$
If $ G $ is non-trivial, then we observe the following:
\begin{enumerate}[font=\normalfont,leftmargin=*]
	\item $ \orbv $ is an injection and $ v $ is an isomorphism.
	\item $ p^N $ is not injective since the numerical classes of two $ G $-equivariant sheaves $ (E, \phi_1) $ and $ (E, \phi_2) $ on $ X $ are mapped to the same element $ \gamma(E) $ in $ N(X) $. $ p^N $ maps into the $ G $-invariant subspace $ N(X)^G $ of $ N(X) $, but may not map onto $ N(X)^G $ since there may be $ G $-invariant sheaves on $ X $ which are not $ G $-linearizable, i.e., which do not lift to $ G $-equivariant sheaves on $ X $.
	\item $ p^R $ is a projection $ (v,(v_{ij})) \mapsto v $ and hence is not injective. $ p^R $ maps onto the $ G $-invariant subspace $ R(X)^G $ of $ R(X) $.
\end{enumerate}
\vspace{8pt}

We first identify Gieseker moduli spaces of stable sheaves in $ \Coh(\cX) $ with Bridgeland moduli spaces of stable objects in $ \D(\cX) $.
\begin{lem}\label{lem_Bridgeland_vs_Gieseker_K3/G}
	Let $ \cX = [K3/G] $ in Situation \ref{situation_finite_group} with a $ G $-invariant ample line bundle $ H $ on $ X $. Choose an element $ x $ in $ N(\cX). $ Let $ y = p^N x $ in $ N(X)^G $ with Mukai vector $ v = (r, c_1, s) $ in $ R(X)^G $ and $ d = \deg(y) $. Suppose $ r > 0, $ and either one of the following conditions is satisfied:
	\begin{enumerate}[font=\normalfont,leftmargin=2em]
		\item $ \gcd(r, d) = 1. $
		\item $ d > 0 $ and $ \gcd(r, d, s) = 1. $
		\item $ d > 0, $ $ y $ is primitive, and $ H $ is $ y $-generic.
		\item $ (r, c_1) $ is primitive and $ H $ is $ y $-generic.
	\end{enumerate}
	Then there exists an interval $ (a,b) $ in $ \RR $ such that for all $ t \in (a,b) $, there is a $ G $-invariant stability condition $ \sigma_t $ in $ \Stab^\dagger(X) $ which is $ y $-generic and
	\begin{equation*}
	M_\cH^{s}(\cX, x) = M_\cH(\cX, x) \cong M_{\wtilde{\sigma_t}}(\D(\cX),x) = M_{\wtilde{\sigma_t}}^{s}(\D(\cX),x).
	\end{equation*}
\end{lem}
\begin{proof}
	Our choice of element $ y $ and polarization $ H $ satisfies the assumption in Proposition \ref{prop_Bridgeland_vs_Gieseker}, which implies there is an integer $ b $ and a real number $ t_0 $ such that for all $ t \geq t_0 $, there is a stability condition 
	$$ \sigma_t = \sigma_{bh,th} \in \Stab^\dagger(X) $$ 
	such that there is an isomorphism
	$$ M_H^{(s)}(X,x) \cong M_{\sigma_t}^{(s)}(\D(X), x) $$ 
	between Gieseker and Bridgeland moduli spaces on $ X $.
	Moreover, these $ \sigma_t $ are $ G $-invariant since $ H $ is $ G $-invariant. Choose any $ t \geq t_0. $ Take a complex $ \cE $ in $ \D(\cX) $ with numerical class $ x $ in $ N(\cX) $. Then we have a chain of equivalences:
	\begin{align*}
	& \cE\ \text{is } \wtilde{\sigma_t}\text{-semistable in } \D(\cX)\\
	\Leftrightarrow \quad & p^*\cE\ \text{is } \sigma_t\text{-semistable in } \D(X) & \text{by Proposition \ref{prop_Bridgeland_relation}}\\
	\Leftrightarrow \quad & p^*\cE \ \text{is } H\text{-semistable in } \Coh(X) & \text{by Proposition \ref{prop_Bridgeland_vs_Gieseker}}\\
	\Leftrightarrow \quad & \cE \ \text{is } \cH\text{-semistable in } \Coh(\cX) & \text{by Proposition \ref{prop_Gieseker_relation}}
	\end{align*}
	Hence there is an isomorphism 
	$$ M_\cH(\cX, x) \cong M_{\wtilde{\sigma_t}}(\D(\cX),x) $$
	between Gieseker and Bridgeland moduli space on $ \cX $.
	
	Since the set of walls determined by $ y $ is locally finite, every interval in $ \RR $ intersects finitely many of them. Therefore, we can choose an interval $ (a,b) $ with $ a \geq t_0 $ such that no walls cross it. Then every $ t \in (a,b) $ determines a stability condition $ \sigma_t $ which is $ y $-generic. 
	
	By our assumptions, there are neither strictly semistable complexes in $ \D(X) $ nor strictly semistable sheaves in $ \Coh(X) $ with numerical class $ y. $ Note that Proposition \ref{prop_Bridgeland_relation} gives an implication
	\begin{equation}\label{implication1}
	p^*\cE\ \text{is } \sigma_t\text{-stable in } \D(X) \ \Rightarrow \ \cE \ \text{is } \wtilde{\sigma_t}\text{-stable in } \D(\cX),
	\end{equation}
	and Proposition \ref{prop_Gieseker_relation} gives another implication
	\begin{equation}\label{implication2}
	p^*\cE\ \text{is } H\text{-stable in } \Coh(X) \ \Rightarrow \ \cE \ \text{is } \cH\text{-stable } \Coh(\cX).
	\end{equation}
	Take a complex $ E $ in $ \D(X) $ with numerical class $ y $ in $ N(X). $ Then Proposition \ref{prop_Bridgeland_vs_Gieseker} also gives an equivalence:
	\begin{equation}\label{equivalence}
	E\ \text{is } \sigma_t\text{-stable in } \D(X) \ \Leftrightarrow \ E \ \text{is } H\text{-stable in } \Coh(X).
	\end{equation}	
	By implication (\ref{implication1}), (\ref{implication2}), and equivalence (\ref{equivalence}), we deduce there are neither strictly semistable complexes in $ \D(\cX) $ nor strictly semistable sheaves in $ \Coh(\cX) $ with numerical class $ x $ as well. This completes the proof.
\end{proof}

Now we consider the minimal resolution $ M \to X/G, $ which gives a derived McKay correspondence
\begin{equation*}
\Phi: \D(\cX) \xrightarrow{\sim} \D(M)
\end{equation*}
Recall that we have a space $ \Stab(\cX) $ of stability conditions on the derived category $ \D(\cX), $ as well as a space $ \Stab(M) $ of stability conditions on $ \D(M). $ The two spaces $ \Stab(\cX) $ and $ \Stab(M) $ are complex manifolds which carry distinguished components $ \Stab^\dagger(\cX) $ and $ \Stab^\dagger(M) $ respectively. The derived equivalence $ \Phi $ induces an isomorphism of complex manifolds
\begin{equation*}
\Phi^S: \Stab(\cX) \xrightarrow{\sim} \Stab(M), \quad (Z, \cP) \mapsto (Z \circ \Phi^{-1}, \Phi \circ \cP),
\end{equation*}
where the slicing $ \Phi \circ \cP $ is defined by $ (\Phi \circ \cP)(\theta) = \Phi(\cP(\theta)) $ for each $ \theta \in \RR. $ We don't know whether $ \Phi^S $ preserves distinguished components. But the following result suffices for our purpose.  
\begin{lem}[{\cite[Proposition 6.1]{beckmann2020equivariant}}]\label{lem_distinguished_FM}
	Let $ \cX = [K3/G] $ in Situation \ref{situation_finite_group} with the minimal resolution $ M \to X/G $ where $ M = \GHilb X $. Let $ \sigma $ be a $ G $-invariant stability condition in the distinguished component $ \Stab^\dagger(X) $. Let $ \wtilde{\sigma} $ denote the induced stability condition on $ \D(\cX) $. Then the derived equivalence $ \Phi: \D(\cX) \xrightarrow{\sim} \D(M) $ gives a stability condition $ \Phi^S(\wtilde{\sigma}) $ in the distinguished component $ \Stab^\dagger(M). $
\end{lem}

The derived equivalence $ \Phi $ preserves Bridgeland stabilities.
\begin{lem}\label{lem_isom_Bridgeland_FM}
	Let $ \tau $ be a stability condition on $ \D(\cX) $. Let $ x $ be an element in $ N(\cX). $ Then there is an isomorphism of Bridgeland moduli stacks
	\begin{equation*}
	M_\tau^{(s)}(\D(\cX),x) \xrightarrow{\simeq} M_{\Phi^S(\tau)}^{(s)}(\D(M),\Phi^N(x)), \quad \cE^{\boldsymbol{\cdot}} \mapsto \Phi(\cE^{\boldsymbol{\cdot}}).
	\end{equation*}
\end{lem}
\begin{proof}
	This is a tautology since the derived equivalence $ \Phi $ maps $ \tau $-(semi)stable objects in $ \D(\cX) $ to $ \Phi^S(\tau) $-(semi)stable objects in $ \D(M) $.
\end{proof}
Now we are ready to prove our main theorem. Let's repeat it here.
\begin{thm}\label{thm_K3/G}
	Let $ \cX = [K3/G] $ in Situation \ref{situation_finite_group}. Let $ x $ be an element in $ N(\cX) $ with $ y = p^N x $ in $ N(X)^G. $ Denote the Mukai vector of $ y $ by $ v = (r,c_1,s) $ in $ R(X)^G $ with $ d = \deg(y). $ Suppose $ r > 0 $ and the pair $ (y,H) $ satisfies the following conditions:
	\begin{enumerate}[font=\normalfont,leftmargin=2em]
		\item $ y $ is primitive and $ H $ is $ y $-generic.
		\item $ d > 0 $ or $ \gcd(r,d) = 1 $.
	\end{enumerate}
	Then $ M(\cX,x) = M^s(\cX,x). $ If $ M(\cX,x) $ is non-empty, then it is an irreducible symplectic manifold of dimension $ n = 2 - \inprod{\orbv(x)^2}_{I\cX} $ deformation equivalent to $ \Hilb^{n/2}(X) $.
\end{thm}
\begin{proof}
	Since $ y $ is primitive and $ H $ is $ y $-generic, by Proposition \ref{prop_K3/G_smoothness_dim}, we have
	$$ M(\cX,x) = M^s(\cX,x). $$ 
	Suppose $ M(\cX, x) $ is not empty. Then it is a smooth projective scheme of dimension $ n = 2 - \inprod{\orbv(x)^2}_{I\cX} $, again by Proposition \ref{prop_K3/G_smoothness_dim}. It remains to show $ M(\cX, x) = M_\cH(\cX, x) $ is deformation equivalent to $ \Hilb^{n/2}(X) $, which will then imply it is irreducible. By assumption, $ r > 0, $ and either $ d > 0 $ or $ \gcd(r,d) = 1 $, so Lemma \ref{lem_Bridgeland_vs_Gieseker_K3/G} applies. It tells there is an interval $ (a,b) $ such that for each $ t \in (a,b), $ there exists a $ G $-invariant stability condition $ \sigma_t $ in $ \Stab^\dagger(X) $ which induces a stability condition $ \tau_t $ on $ \D(\cX) $, and there is an isomorphism
	\begin{equation}\label{eq_first_isom}
	M_\cH(\cX,x) \cong M_{\tau_t}(\D(\cX),x)
	\end{equation}
	between Gieseker and Bridgeland moduli spaces on $ \cX $.
	Take any $ t \in (a,b). $ By Lemma \ref{lem_isom_Bridgeland_FM}, there is an isomorphism
	\begin{equation}\label{eq_second_isom}
	M_{\tau_t}(\D(\cX),x) \cong M_{\Phi^S(\tau_t)}(\D(M),\Phi^N(x))
	\end{equation}
	between Bridgeland moduli spaces on $ \cX $ and $ M $.
	
	Since $ y $ is primitive in $ N(X), $ the element $ x $ is also primitive in $ N(\cX) $: if not, then $ x = kx' $ for some integer $ k > 1 $ or $ k < -1 $, which would then imply $ y = ky' $ with $ y' = p^N x' $, which contradicts the primitivity of $ y $. Therefore, under the isomorphism 
	$$ \Phi^N: N(\cX) \xrightarrow{\sim} N(M), $$ 
	the transformed element $ \Phi^N(x) $ is also primitive in $ N(M) $.
	
	By Lemma \ref{lem_distinguished_FM}, the stability condition $ \Phi^S(\tau_t) $ lies in the distinguished component $ \Stab^\dagger(M) $ because $ \sigma_t $ is in $ \Stab^\dagger(X) $. Now, if we vary $ t $ in $ (a,b), $ the stability conditions $ \sigma_t $ form a real curve in the manifold $ \Stab(X)^G $, which maps to another one in $ \Stab(\cX) $ via the embedding 
	$$ \Stab(X)^G \into \Stab(\cX), \quad \sigma_t \mapsto \tau_t. $$ 
	Under the isomorphism 
	$$ \Phi^S: \Stab(\cX) \xrightarrow{\sim} \Stab(M), $$ 
	the stability conditions $ \Phi^S(\tau_t) $ traces a real curve in the manifold $ \Stab(M) $ as $ t $ varies in $ (a,b). $ Choose some $ c \in (a,b) $ such that $ \Phi^S(\tau_c) $ doesn't lie on any of the walls determined by $ \Phi^N(x). $ Now we have a primitive Mukai vector $ \Phi^N(x) $ in $ N(M) $ and a stability condition $ \Phi^S(\tau_c) $ in the distinguished component $ \Stab^\dagger(M) $ which is $ \Phi^N(x) $-generic. By Theorem \ref{thm_bottini}, the Bridgeland moduli space 
	$$ M_{\Phi^S(\tau_c)}(\D(M),\Phi^N(x)) $$ 
	is deformation equivalent to a Hilbert scheme of points on a $ K3 $ surface, and hence so is the Gieseker moduli space $ M_\cH(\cX,x) $ by the two isomorphisms (\ref{eq_first_isom}) and (\ref{eq_second_isom}).
\end{proof}

If $ G = \{1\} $ is trivial, then Theorem \ref{thm_K3/G} becomes

\begin{cor}\label{cor_K3}
	Let $ (X,H) $ be a polarized $ K3 $ surface. Let $ x $ be an element in $ N(X) $ with Mukai vector $ v = (r,c_1,s) $ in $ R(X)$ and $ d = \deg(x). $ Suppose $ r > 0 $ and the pair $ (x,H) $ satisfies the following conditions:
	\begin{enumerate}[font=\normalfont,leftmargin=2em]
		\item $ x $ is primitive and $ H $ is $ x $-generic.
		\item $ d > 0 $ or $ \gcd(r,d) = 1 $.
	\end{enumerate}
	Then $ M(X,x) = M^s(X,x). $ If $ M(X,x) $ is non-empty, then it is an irreducible symplectic manifold of dimension $ n = 2 - \inprod{v^2}_{X} $ deformation equivalent to $ \Hilb^{n/2}(X) $.
\end{cor}

\begin{rmk}
	\begin{enumerate}[font=\normalfont,leftmargin=*]
		\item Corollary \ref{cor_K3} is a weak version of Theorem \ref{thm_yoshioka} since we require another condition $ d > 0 $ or $ \gcd(r,d) = 1 $, but that's the price we have paid to identify equivariant moduli spaces of stable sheaves in $ \Coh^G(X) $ as Bridgeland moduli spaces of stable complexes in $ \D(X)_G. $
		\item In Theorem \ref{thm_yoshioka}, there is also a result on the non-emptiness of $ M(X,x) $. In our case, we do not yet have a criterion to tell when an equivariant moduli space of stable sheaves on $ K3 $ surfaces is non-empty.
		\item The proof of Theorem \ref{thm_K3/G} shows that if an element $ x $ in $ N(\cX) $ has rank $ r = 1, $ then the results hold with no assumptions on the pair $ (y,H) $. In particular, this applies to $ G $-equivariant Hilbert schemes of points on $ X $.
	\end{enumerate}
\end{rmk}
\vskip 5pt

\section{Appendix: Mukai pairing and HRR formula}\label{review_mukai}

In this appendix we review a few notions in the intersection theory of schemes, including the Mukai vector and the Mukai pairing, and prove the Hirzebruch-Riemann-Roch (HRR) formula. 

Fix an algebraically closed field $ k $. Let $ X $ be a separated scheme of finite type over $ k $. We first recall the Euler characteristic which appears in the left-hand side of the HRR theorem.
\begin{defn}
	If $ X $ is proper, then the Euler characteristic 
	$$ \chi(X, \ \cdot\ ): K(X) \to \ZZ $$
	is defined by
	\begin{equation*}
		\chi(X, E) = \sum_i (-1)^i \dim H^i(X, E)
	\end{equation*}
	for a coherent sheaf $ E $ on $ X $ and extended linearly.
\end{defn}
Recall the HRR theorem for schemes.
\begin{thm}[HRR]
	Let $ X $ be a proper smooth scheme. For all $ x $ in $ K(X), $ we have
	\begin{equation*}
		\chi(X, x) = \int_X \ch(x) \td_X
	\end{equation*}
	in $ \ZZ. $
\end{thm}
\begin{defn}
	If $ X $ is proper, then the Euler pairing
	\begin{equation*}
		\chi: K(X) \times K(X) \to \ZZ
	\end{equation*}
	is defined by
	\begin{equation*}
		\chi(E,F) = \sum_{i} (-1)^i \dim \Ext^i(E,F)
	\end{equation*}
	for two coherent sheaves $ E $ and $ F $ on $ X $ and extended bilinearly.
\end{defn}
If $ X $ is smooth, then $ X $ has the resolution property and $ K(X) \cong K^0(X), $ so we can define an involution on $ K(X) $ via the involution on $ K^0(X) $ which is induced by taking dual vector bundles. 

Let $ X $ be a separated smooth scheme of finite type.
\begin{defn}
	The involution $ (\ \cdot \ )^\vee: K(X) \to K(X) $ is defined by the composition
	\begin{equation*}
		(\ \cdot \ )^\vee: K(X) \xrightarrow{\beta} K^0(X) \xrightarrow{(\ \cdot \ )^\vee} K^0(X) \xrightarrow{\beta^{-1}} K(X)
	\end{equation*}
	where $ \beta $ is the natural isomorphism $ K(X) \xrightarrow{\sim} K^0(X). $
\end{defn}
The following result establishes the relation between the Euler characteristic and the Euler pairing.
\begin{lem}\label{relation_Euler}
	Let $ X $ be a proper smooth scheme. For all $ x $ and $ y $ in $ K(X), $ we have
	\begin{equation*}
		\chi(x,y) = \chi(X, x^\vee y)
	\end{equation*}
	in $ \ZZ. $ In particular, $ \chi(1,x) = \chi(X,x) $ for all $ x $ in $ K(X) $ where $ 1 = [\cO_X] $.
\end{lem}
\begin{proof}
	By the bilinearity of $ \chi, $ it suffices to consider $ x = [V] $ and $ y =[W] $ for vector bundles $ V $ and $ W $ on $ X. $ Then we have
	\begin{align*}
		\chi(x, y) & = \sum_i (-1)^i \dim \Ext^i(V,W) \\
		& = \sum_i (-1)^i \dim H^i(X,V^\vee \otimes W) & \text{because $ V $ is locally free} \\
		& = \chi(X, V^\vee \otimes W) \\
		& = \chi(X, x^\vee y).
	\end{align*}
\end{proof}
Since $ X $ is smooth, the Chow groups $ A(X) = A^*(X) $ form a graded ring with an intersection product. We define an involution on the rational Chow ring $ A(X)_\QQ $ as follows.
\begin{defn}
	The involution $ (\ \cdot \ )^\vee: A(X)_\QQ \to A(X)_\QQ $ is defined by 
	\begin{equation*}
		v^\vee = \sum (-1)^i v_i
	\end{equation*}
	for all $ v = \sum_i v_i $ in $ A(X)_\QQ $ where each $ v_i \in A^i(X)_\QQ. $
\end{defn}
For a unit $ 1 + x = 1 + \sum_{i\geq 1} x_i $ in a graded ring $ R $ with $ x_i \in R^i $ for $ i \geq 1 $ such that $ x^n = 0 $ for some $ n \geq 1, $ we can define its square root by a Taylor series:
\begin{equation*}
	\sqrt{1+x} = \sum_i \binom{1/2}{i} x^i = 1 + \frac{x}{2} - \frac{x^2}{8} + \cdots
\end{equation*}
A quick computation shows the following properties of the involution $ (\ \cdot \ )^\vee $ on $ A(X)_\QQ $.
\begin{lem}
	For all $ v $ and $ w $ in $ A(X)_\QQ $, we have 
	\begin{equation*}
		(vw)^\vee = v^\vee w^\vee \quad \text{and} \quad \sqrt{v^\vee} = \sqrt{v}^\vee
	\end{equation*}
	in $ A(X)_\QQ $ whenever $ \sqrt{v} $ is defined.
\end{lem}
The Todd class of $ X $
\begin{equation*}
	\td_X = 1 + \frac{c_1}{2} + \frac{c_1^2+c_2}{12} + \cdots
\end{equation*}
is a unit in $ A(X)_\QQ $ with $ c_i = c_i(X) = c_i(TX) $ for $ i > 0 $,
and it has a well-defined square root
\begin{equation*}
	\sqrt{\td_X} = 1 + \frac{c_1}{4} + \frac{c_1^2+4c_2}{96} + \cdots
\end{equation*}
in $ A(X)_\QQ^\times. $
\begin{defn}\label{defn_v}
	The Mukai vector map 
	$$ v: K(X) \to A(X)_\QQ $$
	is defined by
	\begin{equation*}
		v(x) = \ch(x) \sqrt{\td_X}
	\end{equation*}
	for all $ x \in K(X). $
\end{defn}
The following properties of the Chern Character, the Mukai vector, and the involutions on $ K(X) $ and $ A(X)_\QQ $ can be easily checked.
\begin{lem}\label{lem_ch_v_inv}
	For all $ x $ and $ y $ in $ K(X), $ we have:
	\begin{enumerate}[font=\normalfont,leftmargin=2em]
		\item $ v(x+y) = v(x) + v(y) $
		\item $ v(xy) = v(x)\ch(y) $
		\item $ \ch(x^\vee) = \ch(x)^\vee $
		\item $ v(x^\vee) = v(x)^\vee \sqrt{\td_X/\td_X^\vee} $
	\end{enumerate}
\end{lem}
\begin{rmk}
	We have the identity
	$$ \sqrt{\td_X/\td_X^\vee} = e^{c_1(X)/2} $$
	in $ A(X)_\QQ $ by applying the splitting principle to $ \td_X $
	as in the proof of Lemma 5.41 in \cite{huybrechts2006fourier}.
\end{rmk}
If $ X $ is proper, then we can define a pairing on $ A(X)_\QQ $.
\begin{defn}\label{defn_v_pairing}
	Let $ X $ be a proper smooth scheme. The \tb{Mukai pairing}
	$$ \inprod{\cdot\ {,}\ \cdot}: A(X)_\QQ \times A(X)_\QQ \to \QQ $$ 
	is defined by
	\begin{equation*}
		\inprod{v,w} = \int_{X} v^\vee w \sqrt{\td_X/\td_X^\vee}
	\end{equation*}
	for all $ v $ and $ w $ in $ A(X)_\QQ. $
\end{defn}
The HRR theorem for schemes implies the HRR formula (\ref{eq_HRR_schemes}) in terms of the Euler pairing and the Mukai pairing. Another version of formula (\ref{eq_HRR_schemes}) is formula (5.5) in \cite{huybrechts2006fourier} where the inputs are complexes of sheaves on the scheme $ X. $
\begin{thm}[HRR Formula]
	Let $ X $ be a proper smooth scheme. For all $ x $ and $ y $ in $ K(X), $ we have
	\begin{equation}\label{eq_HRR_schemes}
		\chi(x, y) = \inprod{v(x),v(y)}
	\end{equation}
	in $ \ZZ. $
\end{thm}
\begin{proof}
	Take any $ x $ and $ y $ in $ K(X). $ We then have
	\begin{align*}
		\chi(x,y) & = \chi(X, x^\vee y) & \text{by Lemma \ref{relation_Euler}} \\
		& = \int_{X} \ch(x^\vee y) \td_X & \text{by the HRR theorem for schemes} \\
		& = \int_{X} \ch(x^\vee) \ch(y) \td_X & \text{because $ \ch $ is a ring map} \\
		& = \int_{X} v(x^\vee) v(y) & \text{by Definition \ref{defn_v}} \\
		& = \int_{X} v(x)^\vee v(y) \sqrt{\td_{X}/\td_{X}^\vee} & \text{by Lemma \ref{lem_ch_v_inv}} \\
		& = \inprod{v(x),v(y)} & \text{by Definition \ref{defn_v_pairing}.}
	\end{align*}
\end{proof}
For a separated smooth scheme $ X $ of finite type, the Chern character map 
$$ \ch: K(X) \to A(X)_\QQ $$ 
is a ring homomorphism and becomes a $ \QQ $-algebra isomorphism after tensored with $ \QQ, $ so the Mukai vector map $ v: K(X) \to A(X)_\QQ $
becomes an isomorphism of $ \QQ $-vector spaces after tensored with $ \QQ $, since $ \sqrt{\td_X} $ is a unit in $ A(\cX)_\QQ. $ If $ X $ is proper, then we can extend the Euler pairing to $ K(X)_\QQ, $ and hence we have two vector spaces $ K(X)_\QQ $ and $ A(\cX)_\QQ $ with bilinear forms. Therefore, we have the following
\begin{prop}
	Let $ X $ be a proper smooth scheme. The Mukai vector map 
	$$ v: K(X) \to A(X)_\QQ $$ 
	induces a linear isometry
	\begin{equation*}
		v: \left(K(X)_\QQ, \chi\right) \xrightarrow{\simeq} \left(A(X)_\QQ, \inprod{\cdot\ {,}\ \cdot}\right).
	\end{equation*}
\end{prop}

\vskip 8pt


\end{document}